\documentclass[12pt]{article}

\usepackage{amsfonts}
\usepackage{graphics}
\usepackage{amsmath}
\usepackage{amscd}
\usepackage{rlepsf}

\newtheorem{Theorem}{Theorem}

\newtheorem{Definition}[Theorem]{Definition}

\newtheorem{Lemma}[Theorem]{Lemma}
\newtheorem{Remark}[Theorem]{Remark}
\newtheorem{Proposition}[Theorem]{Proposition}

\newtheorem{Warning}[Theorem]{Warning}
\newtheorem{Exercise}{Exercise}
\newtheorem{Fundamental Theorem}{Fundamental Theorem}

\newenvironment{Proof}[1][Proof]{\textbf{#1.} }{\ \rule{0.5em}{0.5em}}

\def \Q {\mathbb{Q}}
\def \Z {\mathbb{Z}}

\def \ra {\xrightarrow}

\def \tn {\otimes}
\def \d {\partial}
\def \t {\triangleright}

\def \k {\kappa}

\def \C {\mathcal{C}}

\def \G {\mathcal{G}}

\def \R {\mathbb{R}}

\def \V {\mathcal{V}}

\def \Ga {\Gamma}

\def \f {\phi}

\def \e {\varepsilon}
\def \n {\nu}

\def \s  {\scriptstyle}
\def \S {\Sigma}
\def \R {\mathbb{R}}
\begin{document}

\title{Categorical Groups, Knots and Knotted Surfaces}

\author{Jo\~{a}o  Faria Martins\footnote{email address: jmartins@math.ist.utl.pt}}

\date{\today}
\maketitle
\begin{abstract}
We define a knot invariant and a $2$-knot invariant from any finite
categorical group. We calculate an explicit example for the Spun Trefoil.  
\end{abstract}
\section*{Introduction}
In \cite{Y}, David Yetter defines an invariant of piecewise linear  manifolds
from any  finite categorical group, or, what is the same, from any finite crossed module
$\G=(G,E,\d,\t)$, see \cite{BM} or \cite{BS}. The meaning of Yetter's construction  was elucidated further by Tim Porter in \cite{P}.

\begin{figure}
\centerline{\relabelbox 
\epsfysize 3cm
\epsfbox{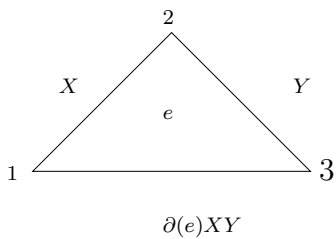}
\relabel {X}{$\s{X}$}
\relabel {Y}{$\s{Y}$}
\relabel {Z}{$\s{\d(e)XY}$}
\relabel {e}{$\s{e}$}
\relabel {1}{$\s{1}$}
\relabel {2}{$\s{2}$}
\relabel {3}{$3$}
\endrelabelbox }
\caption{\label{simplex} Coherent colouring of a simplex. }
\end{figure}

The invariant discovered by Yetter  is defined from triangulations of manifolds equipped with a strict order relation on their set of vertices. Loosely speaking, it  counts, apart from
normalising factors, the number of ways we can colour the edges of the
triangulation by elements of $G$ and  its faces by elements of $E$ in a
coherent way as displayed in figure \ref{simplex}, and so that for each
tetrahedron the faces of each ordered simplex commute as a diagram in $C(\G)$,
the tensor category constructed from $\G$, in other words the categorical group associated with $\G$ (we will explain this construction below).  This last constraint  can be seen
as a  flatness condition. The aim of this article to present a method for calculating this invariant for the case in which the manifold is a knot complement or the $4$-dimensional complement of a knotted surface, therefore defining a class of knot and $2$-knot invariants.  These invariants depend only on the homotopy $2$-type of the complements, similarly to Yetter's invariants of manifolds.

With a bit of work, any manifold invariant defined from triangulations can be
 described in such way  that it can be calculated from handle decompositions of manifolds,
 which  are  more flexible than triangulations\footnote{This is a lesson I learnt from Dr John W. Barrett}. For the case of Yetter's
 invariant, one would consider colourings of $1$-handles by elements of $G$ and
 of $2$-handles by elements of $E$, coherent in a similar fashion to the triangulation scenario above. The
 handles of index $3$ must, like  tetrahedrons, be assigned flatness
 conditions.

 Given a knot diagram or a  movie representation of a knotted surface, there are very
 natural handle decompositions of the complement of any regular neighbourhood  of them.  We will not give
 a precise description of Yetter's invariant in terms of handle
 decompositions of manifolds. Rather,  we will associate a  state sum to every knot diagram
 and any movie of a knotted surface, which is inspired by natural handle decompositions of the
 complements as well as Yetter's construction, and prove directly that they are
 topological invariants.  

 This article is independent of Yetter's results, although it retains  most of
 his construction's philosophy. Nevertheless, it does introduce new techniques and results. 

We include an explicit calculation which asserts that, at least in the $4$-dimensional case, the invariant constructed is non trivial.

It is hoped these descriptions will shed some
 light on how to define invariants of $2$-knots from crossed modules of
 quantum groups. Another possible direction could be trying to define what a
 categorical quandle is.

\section{The Categorical Setting}
\subsection{Crossed Modules and Categorical Groups}
Let $G$ and $E$ be groups. A  crossed module with base $G$ and fibre $E$, say $\G=(G,E,\d,\t)$, is given by a group morphism $\d : E \to G$ and an action $\t$ of $G$ on $E$ on the left by automorphisms. The conditions on $\t$ and $\d$ are:
\begin{enumerate}
\item $\d(X \t e)=X\d(e)X^{-1};\forall X \in G, \forall e \in E$.
\item $\d(e) \t f=e f e^{-1}; \forall e, f \in E$.
\end{enumerate}
Notice that the second condition implies that $\ker \d$ commutes with all $E$.

Let $\G$ be a crossed module. We can define a strict tensor category $\C(\G)$ from it. The set of objects of $\C(\G)$ is given by all elements of $G$. Given a $X \in G$, the set of all morphisms with source $X$ is  $E$, and the target of $e\in E$ is $\d(e) X$. In other words a morphism in $C(\G)$  `` looks like " 
$X \ra{e} \d(e) X$. Given $X \in G$ and $e,f \in E$ the composition

$$X \ra{e} \d(e)X \ra{f} \d(f) \d(e) X$$ is $$X \ra{fe} \d(fe)X.$$   

The tensor product has the form
\begin{equation}\label{Catmon}
\begin{CD}
\begin{CD} X \\ @VVeV \\\d(e)X \end{CD} \tn \begin{CD} Y \\ @VVfV \\\d(f)Y \end{CD} =\begin{CD} XY \\ @VVeX\t f V \\\d(e)X\d(f)Y \end{CD}
\end{CD}.
\end{equation}
From the definition of a crossed module, it is easy to see that we have indeed defined a strict  tensor category. This construction is an old one. The tensor category $\C(\G)$ is a categorical group (see \cite{BM} and \cite{BS}). It is well known that the categories of crossed modules and of categorical  groups are equivalent (see \cite{BS}). We skip further details on this connection since the  point we want to emphasise is that  we can construct a category $\C(\G)$ from any crossed module.   
\subsection{Duality}
The category $\C(\G)$ admits a strict pivotal structure, see \cite{BW}, \cite{BW2} or \cite{FY}, for
which the duality contravariant functor is
\begin{equation}
\left ( X \ra{ e } \d(e) X\right )^*=\left ( X^{-1}\d( e^{-1}) \ra{X^{-1}\t e} X^{-1}\right ).
\end{equation}
From the definition of a crossed module, it is immediate that $*$ is a contravariant functor.
Given an object $X$ of $\C(\G)$, there exists an arrow $\e_X=\left (1_G\ra
  {1_E} X \tn X^*\right )$. With these arrows and the duality contravariant functor $*$, the
  category $\C(\G)$ is a strict pivotal category, in other words verifies the conditions of  definitions $2.1$ and $2.2$ of \cite{BW2}. Equation \ref{Catmon} is used frequently to prove this.  
 Warning: in general this category is not spherical.
\section{A State Sum Invariant of knots}
We now define an invariant of knots, or strictly speaking of dotted oriented
links in which each component has at least one  bivalent vertex. A dotted link is
a link possibly with some extra bivalent vertices inserted.
\subsection{Motivation and Construction}
 Let $D_L \subset \R^2$ be a dotted link
diagram of the dotted link $L$. We thus have a handle decomposition of the 
complement of $L$ for which each arc of the projection generates a $1$-handle, and
each vertex or crossing generates a $2$-handle. We have a unique $0$-handle
around the ''eye of the observer'' and a unique $3$-handle opposite  him/her. Notice that since each link we consider is supposed to have at least a bivalent vertex, there is not the danger that an arc of the projection may be a loop.
This construction is similar to the one in \cite{BGM}, $3.2$. This type of
handle decompositions of knot complements motivates the definition of flat colourings of knot diagrams.

\subsubsection{Colourings of Knot Diagrams}
\begin{figure}
\centerline{\relabelbox
\epsfysize 7cm
\epsfbox{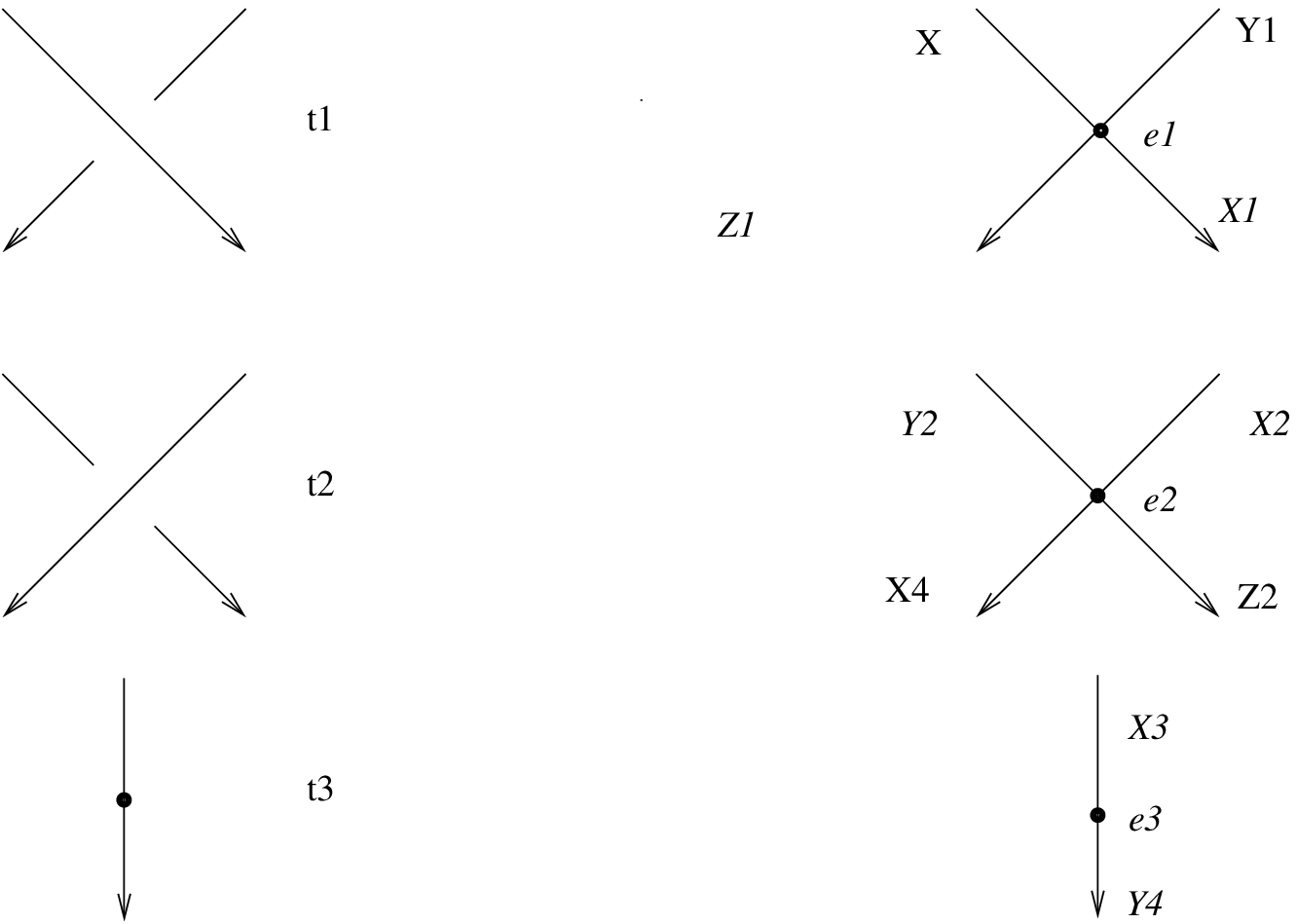}
\relabel {X}{$\s{X}$}
\relabel {X1}{$\s{X}$}
\relabel {X2}{$\s{X}$}
\relabel {X3}{$\s{X}$}
\relabel {X4}{$\s{X}$}
\relabel {Y1}{$\s{Y}$}
\relabel {Y2}{$\s{Y}$}
\relabel {Y4}{$\s{\d (e) X}$}
\relabel {Z1}{$\s{\d(e)XYX^{-1}}$}
\relabel {Z2}{$\s{X^{-1}\d(e)Y X}$}
\relabel {e1}{$\s{e}$}
\relabel {e2}{$\s{e}$}
\relabel {e3}{$\s{e}$}
\relabel {t1}{$\textrm{Yields}$}
\relabel {t2}{$\textrm{Yields}$}
\relabel {t3}{$\textrm{Yields}$}
\endrelabelbox }
\caption{\label{colour1} A bit of a $\G$-colouring of a dotted link diagram.}
\end{figure}

Let $\G=(G,E,\d,\t)$ be a crossed module. Suppose $D_L$ is an oriented dotted link  diagram  of the oriented dotted link $L$. It determines a $4/2$-valent oriented graph $\Gamma(D_L)$ embedded in $\R^2$. 
\begin{Definition}
A $\G$-colouring $c$ of $D_L$ is a colouring $c$ of the edges and the vertices
of the graph $\Gamma(D_L)$ by objects and morphisms of $\C(\G)$, respectively, coherent with
the categorical structure of $\C(G)$, which  in the areas determined
by the  vertices and crossings  looks like figure \ref{colour1}.
\end{Definition} 

\subsubsection{Flat Colourings and the State Sum}\label{flat}
Let $L$ be a dotted link and $D_L$ one of its diagrams. Let also
$\G=(G,E,\d,\t)$ be a crossed module which we assume to be finite. A  $\G$-colouring $c$
of $D_L$, determines a coloured graph $\Gamma(D_L,c)$ embedded in $\R^2$,
such that the edges of $\Ga(D_L)$ are coloured with objects of $\C(\G)$, and such
that at each vertex  of $\Ga(D_L)$ we have a morphism from the tensor product of the objects assigned to
the incoming edges of it  to the tensor product of the objects associated to
the outgoing edges of it. Recall that the category $\C(\G)$ is a strict
pivotal category. Therefore, by a theorem of coherence due to John W. Barrett
and Bruce Westbury, see \cite{BW}, this coloured graph $\Ga(D_L,c)$ can be evaluated to give a
morphism $1_G \ra{\left<\Ga(D_L,c)\right >}1_G$. This morphism only depends on the
planar (but in general not $S^2$) isotopy class of $\Ga(D_L,c)$. See also \cite{FY}.

\begin{Definition}
A colouring $c$ of a dotted link  diagram $D_L$ is said to be flat if:

\begin{equation}
  \left<\Ga(D_L,c)\right >=1_E.
\end{equation}
\end{Definition}
Let $D_L$ be a dotted link diagram of $L$. Define
\begin{equation}
I^3_\G(D_L)=\frac{\#\{ \textrm{flat colourings of } D_L\} }{\#E^{\# \{ \textrm{vertices of } \Ga(D_L)\}}} .\end{equation}
\begin{Exercise}Let $O$ be the unknot with a vertex inserted. Prove that 
$$I^3_\G(O)=\frac{\#G}{\#E}.$$
\end{Exercise}
\begin{Remark} Flatness is not necessary, but ensures non triviality.\end{Remark}

\subsection{Invariance}
We now sketch the proof that $I^3_\G(D_L)$ is a topological invariant.
\subsubsection{Invariance Under Moves of an Isolated Vertex}\label{vertex}
\begin{figure}
\centerline{\relabelbox
\epsfysize 2cm
\epsfbox{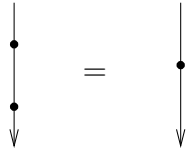}
\relabel{=}{$=$}
\endrelabelbox}
\caption{\label{vabs} Simplest vertex move.}
\end{figure}

First of all the following holds:
\begin{Lemma}
The state sum $I^3_\G(D_L)$ is invariant under the move of figure \ref{vabs}.
\end{Lemma}
The method of proof of the invariance under this move generalises to the other moves, which is  why we will present it in full detail. Given a dotted link
diagram $D$, let $C(D)$ be the set of flat colourings of $D$.

\begin{figure}
\centerline{\relabelbox
\epsfysize 2cm
\epsfbox{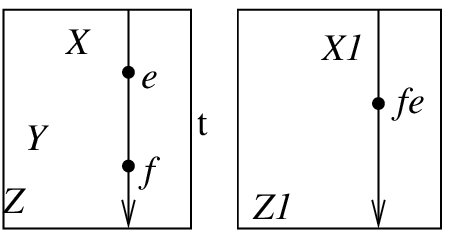}
\relabel{X}{$\s{X}$}
\relabel{X1}{$\s{X}$}
\relabel{Y}{$\s{\d(e)X}$}
\relabel{Z}{$\s{\d(fe)X}$}
\relabel{Z1}{$\s{\d(fe)X}$}
\relabel{e}{$\s{e}$}
\relabel{f}{$\s{f}$}
\relabel{fe}{$\s{fe}$}
\relabel{t}{$\substack{F\\ \mapsto}$}
\endrelabelbox}
\caption{\label{vabsproof} A map between colourings.}
\end{figure}

\begin{Proof}
Let the diagrams $D$ and $D'$ differ by the move in figure
\ref{vabs}. Consider the map $F:C(D)\to C(D')$ defined in figure
\ref{vabsproof}, with no changes outside the rectangle. It is obvious that it
sends flat colourings to flat colourings since the morphisms of $\C(\G)$ inside the
rectangle in $D$ and $D'$ are the same, namely $X \ra{fe} \d(fe) X$. Obviously the map $F$ is surjective
and its fibre at each point (in other words the inverse image at each point) has the same cardinality as $E$, which finishes the proof.  \end{Proof}

\begin{figure}
\centerline{\relabelbox
\epsfysize 4cm
\epsfbox{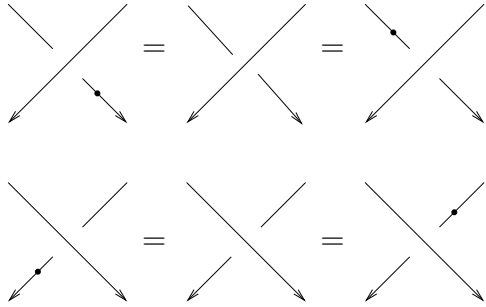}
\relabel{=1}{$=$}
\relabel{=2}{$=$}
\relabel{=3}{$=$}
\relabel{=4}{$=$}
\endrelabelbox}
\caption{\label{v1}Second simplest vertex moves.}
\end{figure}

Secondly:
\begin{Lemma} The state sum $I^3_\G(D)$ is invariant under the moves of figure
  \ref{v1}.\end{Lemma}

\begin{figure}
\centerline{\relabelbox
\epsfysize 4cm
\epsfbox{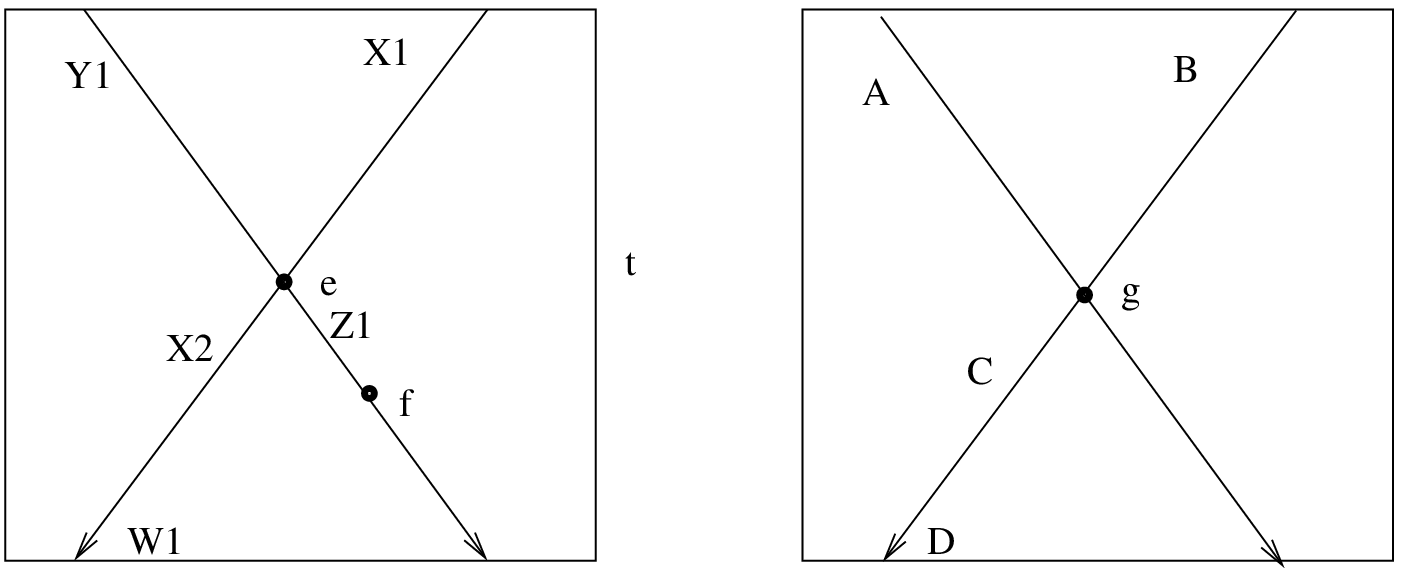}
\relabel{Y1}{$\s{Y}$}
\relabel{X1}{$\s{X}$}
\relabel{X2}{$\s{X}$}
\relabel{Z1}{$\s{X^{-1}\d(e)YX}$}
\relabel{W1}{$\s{\d(f)X^{-1}\d(e)YX}$}
\relabel{A}{$\s{Y}$}
\relabel{B}{$\s{X}$}
\relabel{C}{$\s{X}$}
\relabel{D}{$\s{\d(f)X^{-1}\d(e)YX}$}
\relabel{e}{$\s{e}$}
\relabel{f}{$\s{f}$}
\relabel{g}{$\s{(X \t f) e }$}
\relabel{t}{$\substack{F\\ \mapsto}$}
\endrelabelbox}
\caption{\label{v1proof} A map between colourings.}
\end{figure}

\begin{Proof} Let us prove the upper left  corner. The proof for the other
  cases is analogous.
 Let $D$ and $D'$ be link diagrams differing by the upper left
move in figure \ref{v1}. We define a map $F:C(D)\to C(D')$ as in figure
\ref{v1proof}, where the rest of the colouring remains unaltered.

As before, $F$ transforms flat colourings into flat colourings. This is because
in both cases the morphism in $\C(\G)$ inside the square is the same because of
equation (\ref{Catmon}). This morphism is $YX \ra{(X \t f) e} \d(X\t f)\d(e)
XY =X \d(f)X^{-1} \d(e) Y X$. Obviously $F$ is surjective and the fibre at each
point has the same cardinality as $E$.
\end{Proof}
\begin{figure}
\centerline{\relabelbox
\epsfysize 2cm
\epsfbox{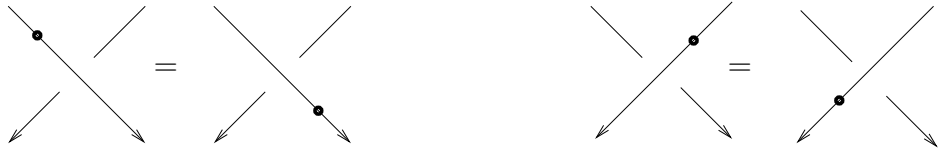}
\relabel{=1}{$=$}
\relabel{ =2}{$=$}
\endrelabelbox}
\caption{\label{v2} Remaining vertex moves.}
\end{figure}

Next we have:
\begin{Lemma}
The state sum $I^3_\G(D)$ is invariant under the moves of figure \ref{v2}.
\end{Lemma}
\begin{Proof}
We prove the equation on the left. Let $D$ and $D'$ differ by this move. We define a map
$F:C(D) \to C(D')$ as in figure \ref{v2proof}.
\begin{figure}
\centerline{\relabelbox
\epsfysize 4cm
\epsfbox{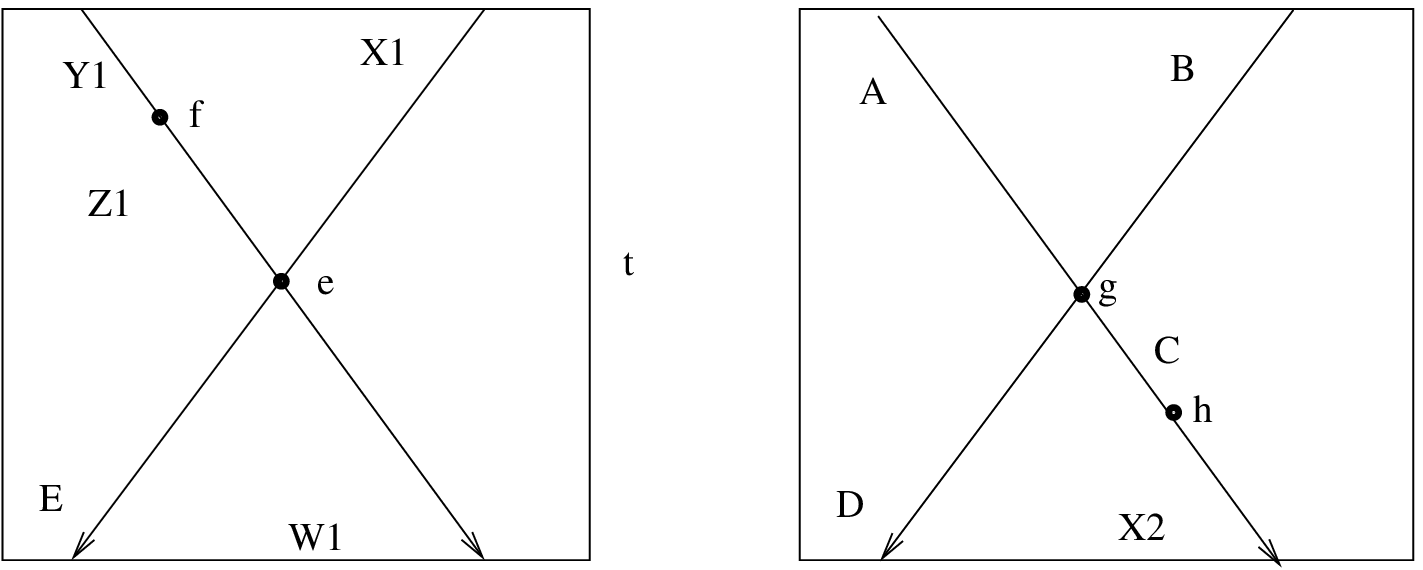}
\relabel{Y1}{$\s{Y}$}
\relabel{X1}{$\s{Z}$}
\relabel{Z1}{$\s{\d(e)X}$}
\relabel{E}{$\s{W}$}
\relabel{W1}{$\s{\d(e)X}$}
\relabel{A}{$\s{X}$}
\relabel{B}{$\s{Z}$}
\relabel{C}{$\s{X}$}
\relabel{D}{$\s{W}$}
\relabel{X2}{$\s{\d(e)X}$}
\relabel{e}{$\s{f}$}
\relabel{f}{$\s{e}$}
\relabel{g}{$\s{feXZX^{-1}\t e^{-1}}$}
\relabel{h}{$\s{ e }$}
\relabel{t}{$\substack{F\\ \mapsto}$}
\endrelabelbox}
\caption{\label{v2proof} A map between colourings. Here $W=\d(f)\d(e)XZX^{-1}\d(e^{-1})=\d(f)\d(e)\d(XZX^{-1}\t e^{-1})XZX^{-1}$.}
\end{figure}
To prove it transforms flat colourings into flat colourings notice that $(W\t e) feXZX^{-1}\t e^{-1}=fe$ (c.f. equation (\ref{Catmon})). Obviously $F$ is a bijection.
\end{Proof}

\subsubsection{Invariance Under Reidemeister-I}
\begin{figure}
\centerline{\relabelbox 
\epsfysize 3cm
\epsfbox{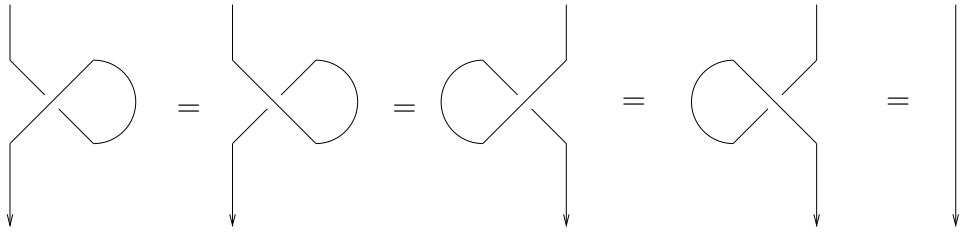} 
\relabel{=1}{$=$}
\relabel{=2}{$=$}
\relabel{=3}{$=$}
\relabel{=4}{$=$}
\endrelabelbox }
\caption{Different Cases of Reidemeister-I.}
\label{R1}
\end{figure}

Since we are considering oriented knot diagrams, there are four different cases of the  Reidemeister-I move, considered up to planar isotopy. They are depicted in figure \ref{R1}. 

This Reidemeister move is the easiest to verify. This is because of the
identities of figure \ref{R1}.  Therefore, we can apply the results of
\ref{vertex}. In particular they tell us that vertices can be absorbed by the
rest of the diagram. We omit the details since more intricate calculations
will follow. 
\begin{figure}
\centerline{\relabelbox 
\epsfysize 3cm
\epsfbox{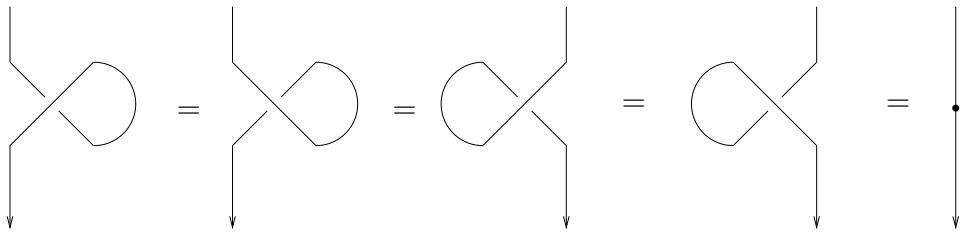} 
\relabel{=1}{$=$}
\relabel{=2}{$=$}
\relabel{=3}{$=$}
\relabel{=4}{$=$}
\endrelabelbox }
\caption{Identity used to prove invariance under Reidemeister-I.}
\label{r1}
\end{figure}

\subsubsection{Invariance Under Reidemeister-II}
\begin{figure}
\centerline{\relabelbox 
\epsfysize 3cm
\epsfbox{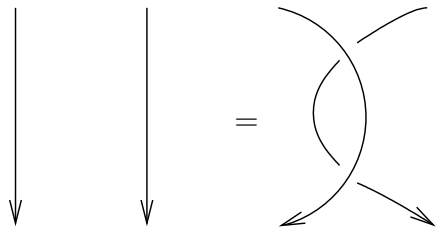} 
\relabel{=}{$=$}
\endrelabelbox }
\caption{One variant of the Reidemeister-II move.}
\label{R2}
\end{figure}

There are four different kinds of oriented Reidemeister-II moves. They are
obtained from figure \ref{R2} through considering all the possible
orientations of  the two strands.

\begin{figure}
\centerline{\relabelbox 
\epsfysize 3cm
\epsfbox{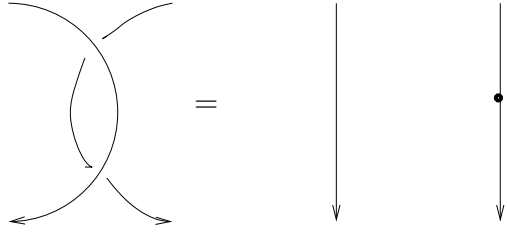} 
\relabel{=}{$=$}
\endrelabelbox }
\caption{Identity used to prove invariance under Reidemeister-II move.}
\label{R2p}
\end{figure}

\begin{figure}
\centerline{\relabelbox 
\epsfysize 3cm
\epsfbox{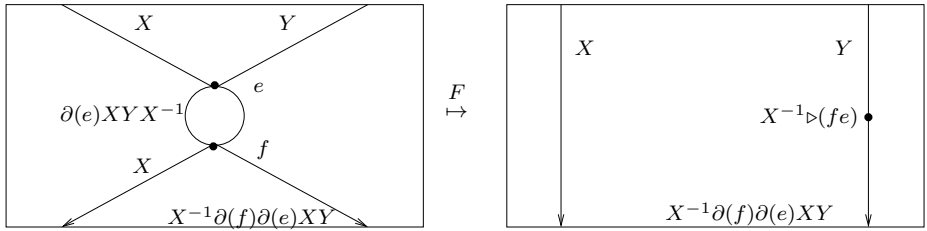} 
\relabel{t}{$\substack{F \\ \mapsto}$}
\relabel{A}{$\s{X}$}
\relabel{B}{$\s{Y}$}
\relabel{C}{$\s{X}$}
\relabel{D}{$\s{\d(e)XYX^{-1}}$}
\relabel{E}{$\s{X^{-1}\d(f)\d(e)XY}$}
\relabel{F}{$\s{X}$}
\relabel{G}{$\s{Y}$}
\relabel{H}{$\s{X^{-1}\d(f)\d(e)XY}$}
\relabel{e}{$\s{e}$}
\relabel{f}{$\s{f}$}
\relabel{g}{$\s{X^{-1} \t (fe)}$}
\endrelabelbox }
\caption{Map used to prove invariance under Reidemeister-II.}
\label{R2proof}
\end{figure}

Let $D$ and $D'$ be knot diagrams differing by a Reidemeister move of type
II. We can suppose by \ref{vertex} that both diagrams have at least one
vertex. Therefore, all is obvious from the identity in figure \ref{R2p} (and its counterparts
for other orientations of the strands), together with \ref{vertex}. To prove it use the
map $F:C(D) \to C(D')$ of figure \ref{R2proof}. Note that it sends flat
colourings to flat colourings since in both cases the morphism inside the
rectangle  is $XY \ra{fe} \d(fe) XY$. As before $F$ is surjective and its fibre has the same cardinality as $E$.
\subsubsection{Invariance Under Reidemeister-III}
\begin{figure}
\centerline{\relabelbox 
\epsfysize 4cm
\epsfbox{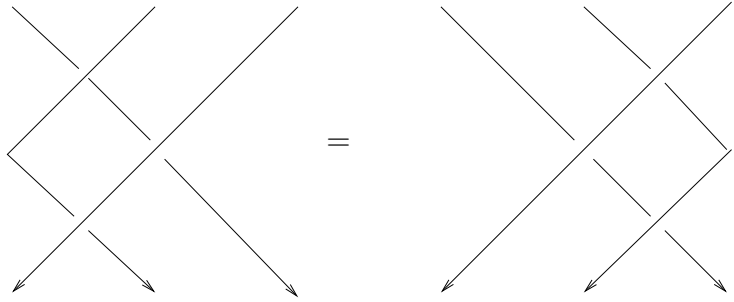}
\relabel{t}{$=$}
\endrelabelbox}
\caption{Reidemeister-III.}
\label{R3}
\end{figure}

This is the most difficult move. Since we are considering oriented knot diagrams, there are sixteen versions of
the Reidemeister move number III. They are obtained from the move in figure \ref{R3},
through considering the mirror image as well as all possible orientations of
the strands.

\begin{figure}
\centerline{\relabelbox 
\epsfysize 4cm
\epsfbox{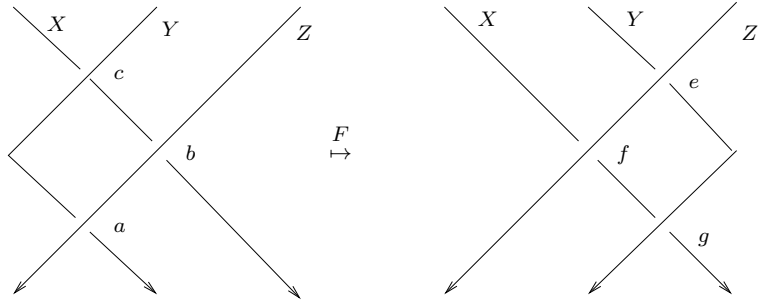}
\relabel{t}{$\substack{F \\ \mapsto}$}
\relabel{X}{$\s{X}$}
\relabel{X1}{$\s{X}$}
\relabel{Y}{$\s{Y}$}
\relabel{Y1}{$\s{Y}$}
\relabel{Z}{$\s{Z}$}
\relabel{Z1}{$\s{Z}$}
\relabel{e}{$\s{e}$}
\relabel{f}{$\s{f}$}
\relabel{g}{$\s{g}$}
\relabel{a}{$\s{a}$}
\relabel{b}{$\s{b}$}
\relabel{c}{$\s{c}$}
\endrelabelbox}
\caption{Map used to prove invariance under Reidemeister-III move, first case.}
\label{R3proof}
\end{figure}

\begin{figure}
\centerline{\relabelbox 
\epsfysize 4cm
\epsfbox{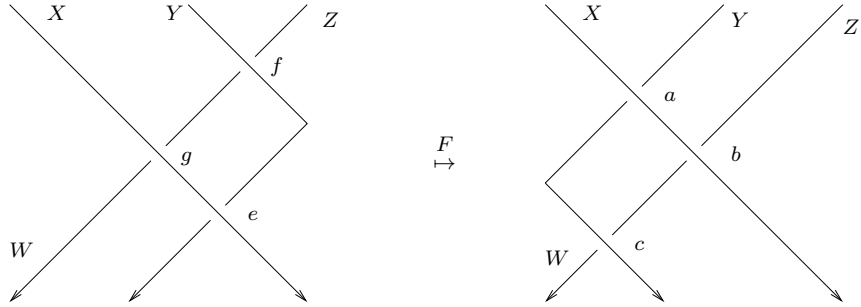}
\relabel{t}{$\substack{F \\ \mapsto}$}
\relabel{X0}{$\s{X}$}
\relabel{X1}{$\s{X}$}
\relabel{Y0}{$\s{Y}$}
\relabel{Y1}{$\s{Y}$}
\relabel{Z0}{$\s{Z}$}
\relabel{Z1}{$\s{Z}$}
\relabel{e}{$\s{f}$}
\relabel{f}{$\s{e}$}
\relabel{g}{$\s{g}$}
\relabel{a}{$\s{a}$}
\relabel{b}{$\s{b}$}
\relabel{c}{$\s{c}$}
\relabel{W1}{$\s{W}$}
\relabel{W2}{$\s{W}$}
\endrelabelbox}
\caption{Map used to prove invariance under Reidemeister-III move, second case. Here $W=\d(c)\d(a)XYX^{-1} \d(b)XZY^{-1}X^{-1}\d(a^{-1})$.}
\label{MapR32}
\end{figure}

 Let $D$ and $D'$ differ by a Reidemeister-III. Suppose it is of the kind
 depicted in figure \ref{R3}. We leave it as an exercise to prove that the the equalities $a=e$, $b=(Y^{-1}\t e^{-1})(Y^{-1} Z \t g)$ and $c= f ( X\t e)$ 
 in figure \ref{R3proof}, define a one-to-one map between colourings of $D$ and
 colourings of $D'$. Note that we always have $ a (Y \t b) c= (Z \t g )f(X
 \t  e)$ (recall equation (\ref{Catmon})), which implies that flat colourings are always sent to flat
 colourings under the map just defined. The proof of invariance for the other
 types of Reidemeister-III is analogous. For example for the mirror image of
 figure \ref{R3}, we use the map of figure \ref{MapR32} to prove
 invariance. Here $g=ca$, $f=YX^{-1} \t b YZY^{-1} X^{-1} \t a ^{-1}$ and
 $e=a$. This map transforms flat colourings into flat colourings since $c
 (\d(a)XYX^{-1}\t b) a=ca XYX^{-1} \t b=(W \t e) g (X \t f)$.

\subsubsection{Planar Deformations}
 The invariance of $I^3_\G(D)$ under horizontal deformations is an obvious
 consequence of the general theorem of coherence for pivotal categories, see
 \cite{BW} or \cite{FY}. As mentioned in \ref{flat}, evaluations of colourings
 are invariant under planar deformations, thus flat colourings are sent to
 flat colourings under planar deformations. This finishes the proof that $I^3_\G(D)$ is a  topological invariant. We thus have a knot invariant $I^3_\G$ for any finite crossed module $\G=(G,E,\d,\t)$. The discussion of it appears in section $4$.
\begin{Remark}
Despite the fact that we need to consider oriented knot diagrams to define our invariant, it is easy to verify that our invariant does not depend on the orientations chosen.
\end{Remark}

\section{Knotted Surfaces}
As usual, we fix a finite crossed module $\G=(G,E,\d,\t)$.

For details on knotted surfaces see \cite{CS}. Consider an embedded surface
$\S$ in $S^4=\R^4 \cup \{\infty\}$, or in  general an embedding in codimension
$2$. Suppose the projection in the last variable is a Morse function in $\S$,
thus determining a handle decomposition of $\S$. Let $\n(\S)$ be a  regular
neighbourhood of $\S$ in $S^4$. We then have a natural handle
decomposition of $S^4 \setminus \n(\S)$, where a regular neighbourhood of each
$a$-handle of $\S$ intersects complementary a unique $a+1$ handle of $S^4 \setminus \n(\S)$, see \cite{GS}, $6.2$. This is very easy to visualise in
dimension $3$. Therefore, according to the discussion in the introduction, if we consider the movie of the chosen embedding
of $\S$, births of circles will induce $1$-handles in the handle decomposition  of the
complement, and therefore must be assigned elements of $G$; saddle points
induce $2$-handles of the complement, thus must be coloured by elements of
$E$, coherently; and deaths of circles correspond to $3$-handles in the
complement, thus must correspond to flatness conditions.  

However, in general when considering the Kirby Diagram of this
handle decomposition, the attaching regions of  $2$-handles are only
determined up to handle slides, which causes an extra complication in our
combinatorial framework. This problem is solved below by considering a set of
relations on dotted knot diagrams.

\subsection{A Vector Space Associated with Knot Diagrams} 
\begin{figure}
\centerline{\relabelbox 
\epsfysize 3cm
\epsfbox{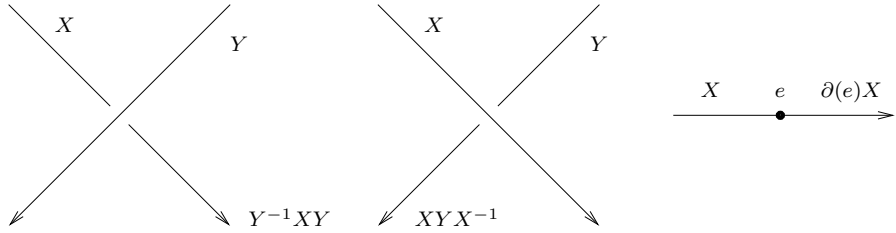} 
\relabel {X1}{$\s{X}$}
\relabel {X2}{$\s{X}$}
\relabel {X}{$\s{X}$}
\relabel {Y1}{$\s{Y}$}
\relabel {Y}{$\s{Y}$}
\relabel {Z}{$\s{Y^{-1}XY}$}
\relabel {W}{$\s{XYX^{-1}}$}
\relabel {e}{$\s{e}$}
\relabel {T}{$\s{\d(e)X}$}
\endrelabelbox }
\caption{Definition of a colouring of a dotted knot diagram.}
\label{Colour}
 \end{figure}

 Recall that a  dotted knot diagram is
 by definition a regular projection of a bivalent graph, in other words of a
 link possibly with some extra bivalent vertices inserted. 
Let $D$ be a dotted knot diagram as before oriented. 
 Let also  $\G=(G,E,\t,\d)$ be a  finite crossed module.
\begin{Definition}
A colouring of $D$ is an assignment of elements of $G$ to the arcs of $D$ and elements of $E$ to the vertices of $D$ verifying the conditions of figure \ref{Colour}.\end{Definition}
\begin{Warning} The definition of a colouring of a dotted knot diagram differs from  the one given previously in the $3$-dimensional case.
\end{Warning}
\begin{Remark}\label{refer}
As in the $3$-dimensional case, a $\G$-colouring of a dotted knot diagram $D$
determines a colouring of the vertices and the edges of the $4/2$-valent graph
determined by $D$ by objects and morphisms of $\C(G)$, the tensor category
made from $\G$. However, in this case, the $4$-valent vertices are always
assigned the trivial morphism. Nevertheless, the important observation is
that, as in  \ref{flat}, colourings of
dotted knot diagrams can always be evaluated to give morphisms in $\C(\G)$.
\end{Remark}
\begin{Definition}
Let $D$ be a knot diagram. A dotting of $D$ is an insertion of bivalent vertices on $D$ considered up to planar isotopy. If $D$ is a knot diagram, let $V(D)$ be the free $\Q$-vector space on the set of all colourings of all dottings of $D$.  \end{Definition}

\begin{figure}
\centerline{\relabelbox 
\epsfysize 3cm
\epsfbox{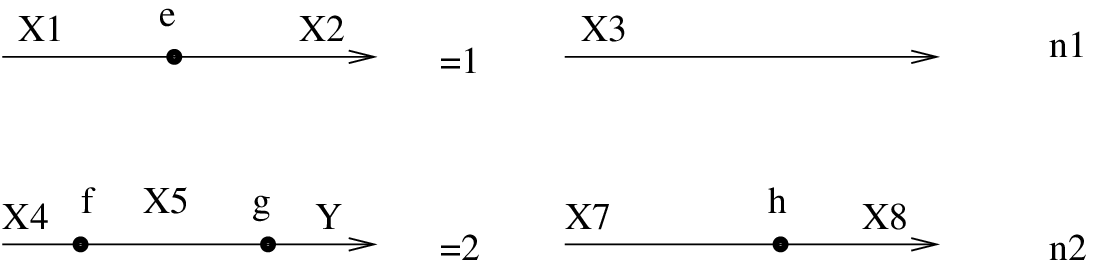} 
\relabel {X1}{$\s{X}$}
\relabel {X2}{$\s{X}$}
\relabel {X3}{$\s{X}$}
\relabel {X4}{$\s{X}$}
\relabel {X5}{$\s{\d(e)X}$}
\relabel {Y}{$\s{\d(fe)X}$}
\relabel {X7}{$\s{X}$}
\relabel {X8}{$\s{\d(fe)X}$}
\relabel {e}{$\s{1_E}$}
\relabel {f}{$\s{e}$}
\relabel {g}{$\s{f}$}
\relabel {h}{$\s{fe}$}
\relabel{=2}{$=$}
\relabel{=1}{$=$}
\relabel{n1}{$R1$}
\relabel{n2}{$R2$}
\endrelabelbox }
\caption{Relations on colourings.}
\label{Relations1}
\end{figure}
\begin{figure}
\centerline{\relabelbox 
\epsfysize 8cm
\epsfbox{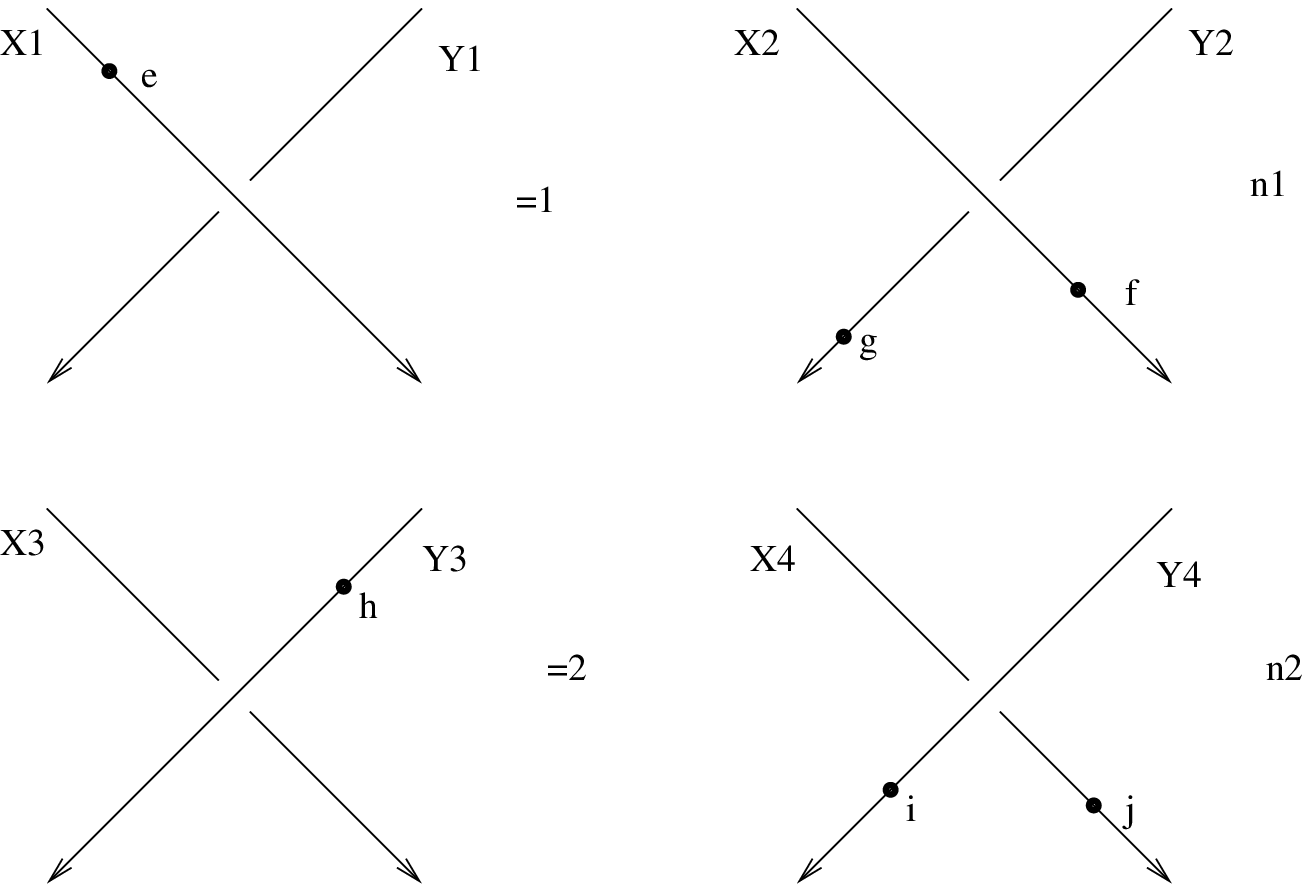} 
\relabel {X1}{$\s{X}$}
\relabel {X2}{$\s{X}$}
\relabel {X3}{$\s{X}$}
\relabel {X4}{$\s{X}$}
\relabel {Y1}{$\s{Y}$}
\relabel {Y2}{$\s{Y}$}
\relabel {Y3}{$\s{Y}$}
\relabel {Y4}{$\s{Y}$}
\relabel {e}{$\s{e}$}
\relabel {f}{$\s{e}$}
\relabel {g}{$\s{eXYX^{-1}\t e^{-1}}$}
\relabel {h}{$\s{e}$}
\relabel {i}{$\s{e}$}
\relabel {j}{$\s{Y^{-1}\t e^{-1}Y^{-1}X \t e}$}
\relabel{=2}{$=$}
\relabel{=1}{$=$}
\relabel{n1}{$R3$}
\relabel{n2}{$R4$}
\endrelabelbox }
\caption{Relations on colourings.}
\label{Relations2}
\end{figure}
\begin{figure}
\centerline{\relabelbox 
\epsfysize 8cm
\epsfbox{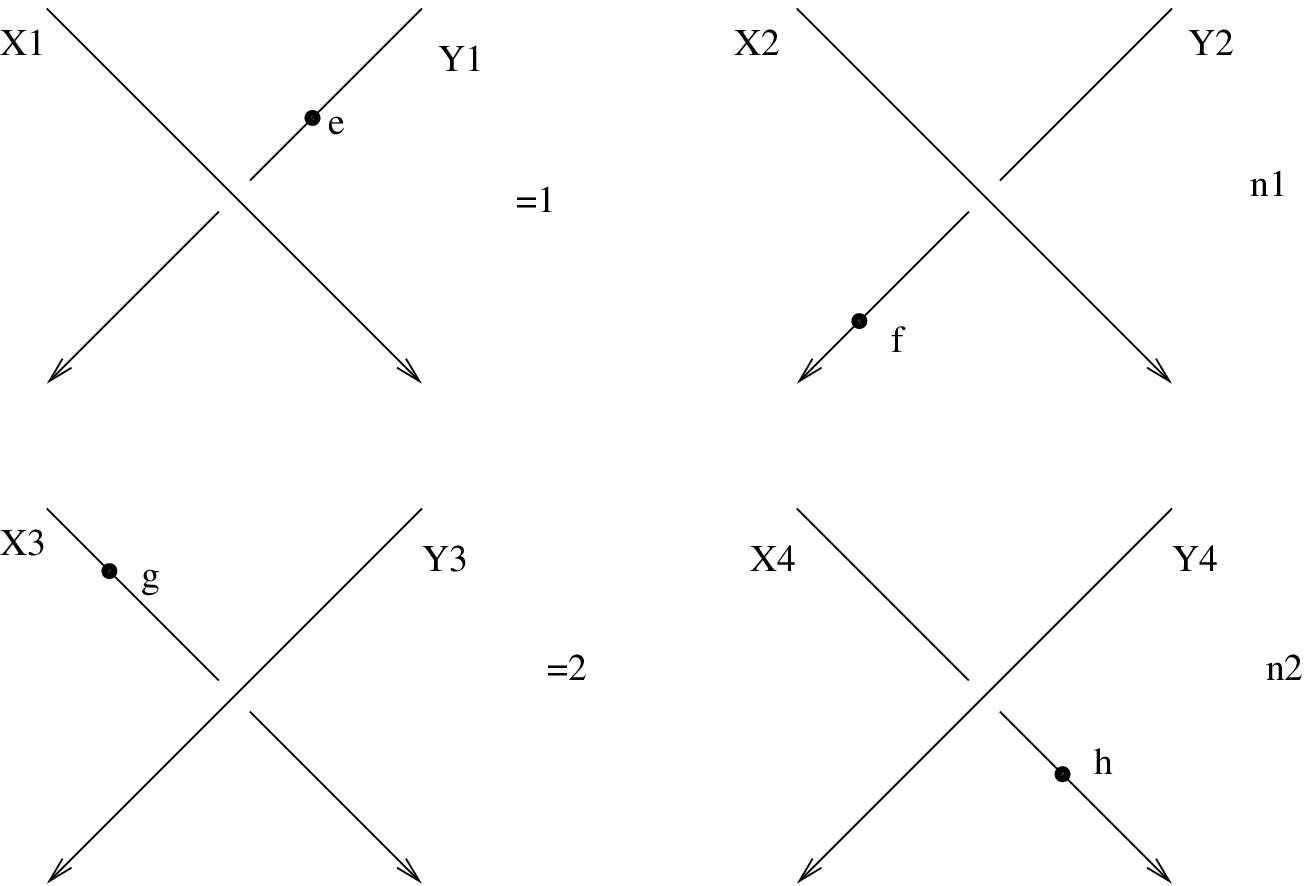} 
\relabel {X1}{$\s{X}$}
\relabel {X2}{$\s{X}$}
\relabel {X3}{$\s{X}$}
\relabel {X4}{$\s{X}$}
\relabel {Y1}{$\s{Y}$}
\relabel {Y2}{$\s{Y}$}
\relabel {Y3}{$\s{Y}$}
\relabel {Y4}{$\s{Y}$}
\relabel {e}{$\s{e}$}
\relabel {f}{$\s{X\t e}$}
\relabel {g}{$\s{e}$}
\relabel {h}{$\s{Y^{-1} \t e}$}
\relabel{=2}{$=$}
\relabel{=1}{$=$}
\relabel{n1}{$R5$}
\relabel{n2}{$R6$}
\endrelabelbox }
\caption{Relations on colourings.}
\label{Relations3}
\end{figure}

Consider now the relations of figures \ref{Relations1}, \ref{Relations2} and
\ref{Relations3}. It is straightforward to see that they are local on the knot
diagrams and that they transform colourings into colourings.
\begin{Remark}\label{argument}
 Notice that in all cases
the morphisms in $\C(\G)$ evaluated considering the colouring of both sides of
the  relations $R1$ to $R6$ (c.f. remark \ref{refer}) are always the same. 
\end{Remark}
\begin{Definition}
Let $D$ be an oriented knot diagram (without vertices). The vector space $\V(D)$ is defined as the vector space obtained from $V(D)$ by considering the relations $R1$ to $R6$.
\end{Definition}

\subsubsection{Some Commutation Relations}\label{commutation}
\begin{figure}
\centerline{\relabelbox 
\epsfysize 8cm
\epsfbox{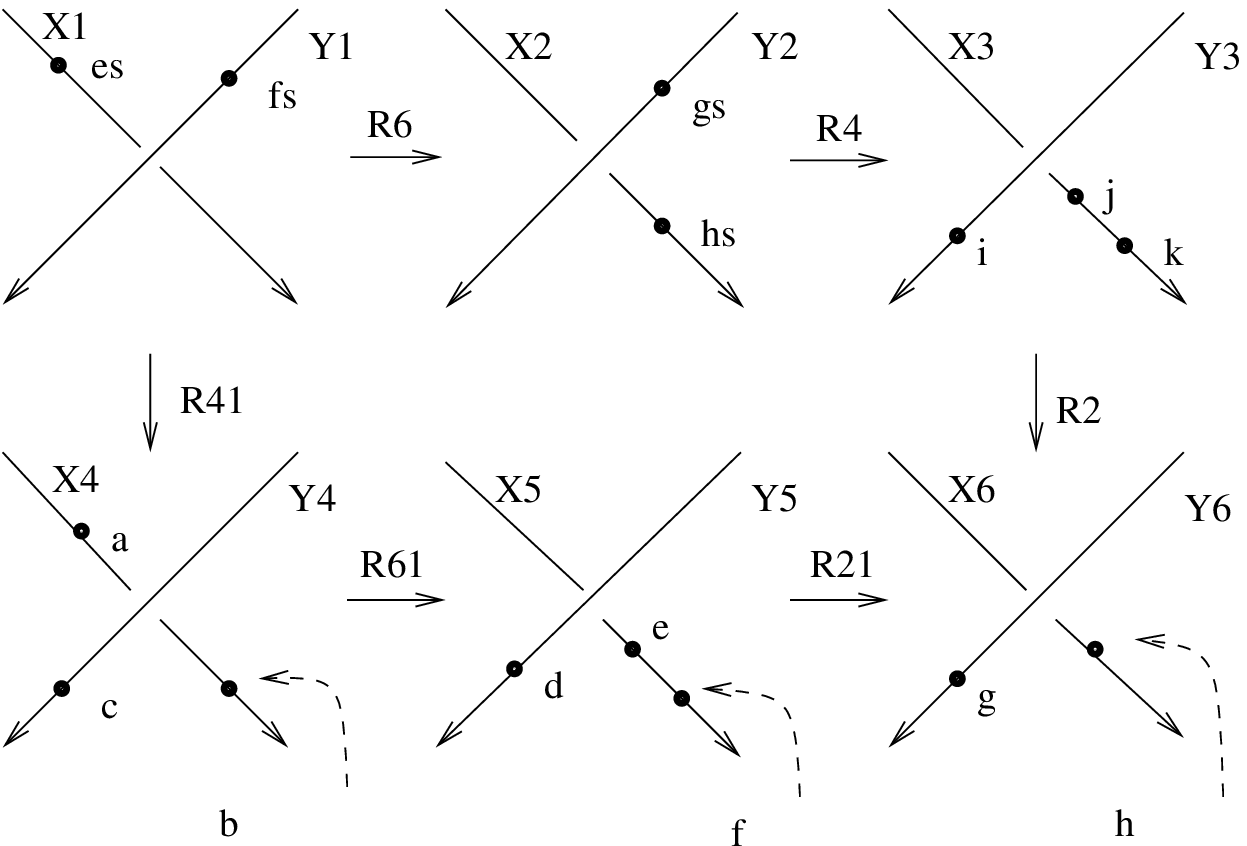}
\relabel{es}{$\s{e}$}
\relabel{fs}{$\s{f}$}
\relabel{gs}{$\s{f}$}
\relabel{hs}{$\s{Y^{-1}\d(f^{-1})\t e}$}
\relabel{i}{$\s{f}$}
\relabel{j}{$\s{Y^{-1}\t f^{-1} Y^{-1}X \t f}$}
\relabel{k}{$\s{Y^{-1}\d(f^{-1})\t e}$}
\relabel{X1}{$\s{X}$}
\relabel{X2}{$\s{X}$}
\relabel{X3}{$\s{X}$}
\relabel{X4}{$\s{X}$}
\relabel{X5}{$\s{X}$}
\relabel{X6}{$\s{X}$}
\relabel{Y1}{$\s{Y}$}
\relabel{Y2}{$\s{Y}$}
\relabel{Y3}{$\s{Y}$}
\relabel{Y4}{$\s{Y}$}
\relabel{Y5}{$\s{Y}$}
\relabel{Y6}{$\s{Y}$}
\relabel{a}{$\s{e}$}
\relabel{b}{$\s{Y^{-1}\t f^{-1} Y^{-1} \d(e)X \t f}$}
\relabel{c}{$\s{f}$}
\relabel{d}{$\s{f}$}
\relabel{e}{$\s{Y^{-1}\t e}$}
\relabel{f}{$\s{Y^{-1}\t f^{-1} Y^{-1}\d(e)X \t f }$}
\relabel{g}{$\s{f}$}
\relabel{h}{$\s{Y^{-1} \t f^{-1} Y^{-1}\t e Y^{-1}X\t f}$}
\relabel{R6}{${R6}$}
\relabel{R4}{${R4}$}
\relabel{R41}{${R4}$}
\relabel{R61}{${R6}$}
\relabel{R21}{${R2}$}
\relabel{R2}{${R2}$}
\endrelabelbox}
\caption{A commutation relation.}
\label{com1}
\end{figure}

\begin{figure}
\centerline{\relabelbox 
\epsfysize 8cm
\epsfbox{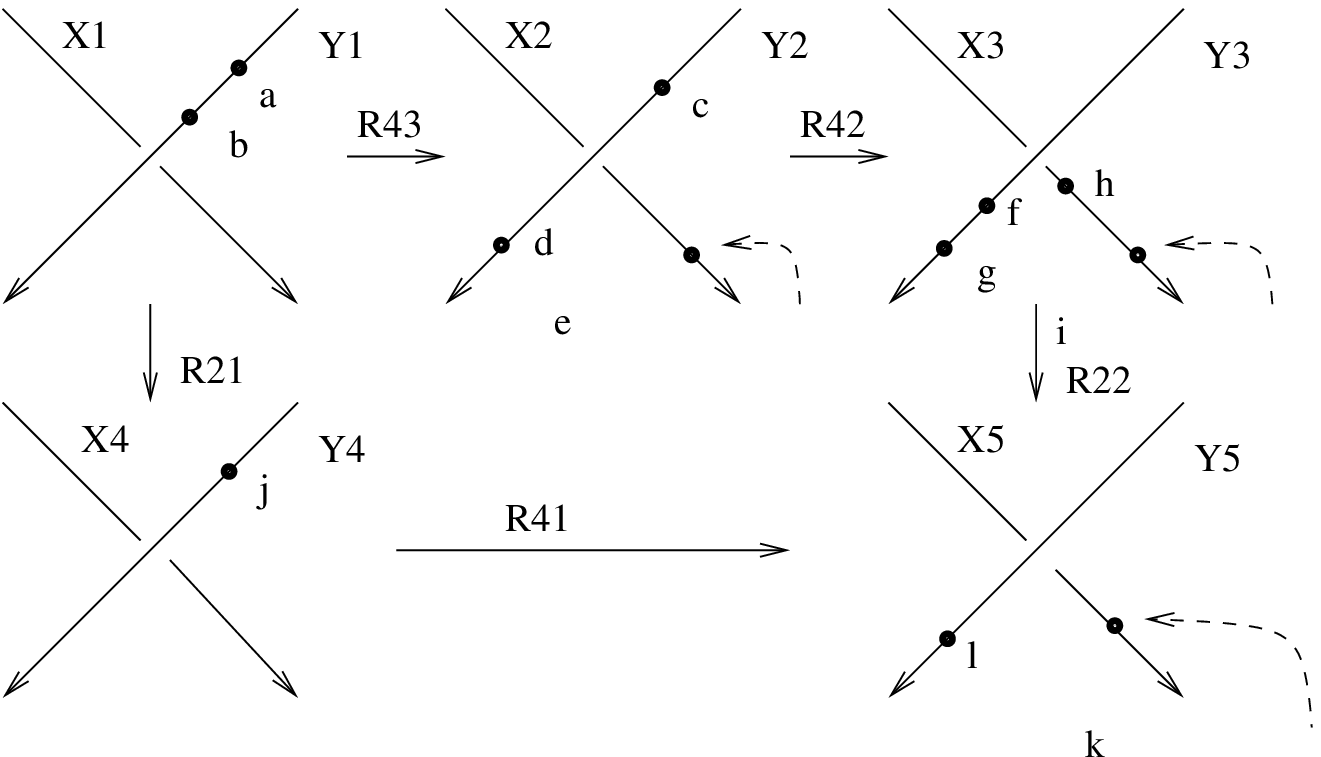}
\relabel{a}{$\s{a}$}
\relabel{b}{$\s{b}$}
\relabel{c}{$\s{a}$}
\relabel{d}{$\s{ b}$}
\relabel{e}{$\s{Y^{-1}\d(a^{-1})\t b^{-1}Y^{-1}\d(a^{-1})X \t b }$}
\relabel{f}{$\s{a}$}
\relabel{g}{$\s{b}$}
\relabel{h}{$\s{Y^{-1}\t a^{-1} Y^{-1}X \t a }$}
\relabel{i}{$\s{Y^{-1}\d(a^{-1})\t b^{-1}Y^{-1}\d(a^{-1})X \t b }$}
\relabel{j}{$\s{ba}$}
\relabel{k}{$\s{Y^{-1}\t (a^{-1} b^{-1}) Y^{-1}X\t (ba)}$}
\relabel{l}{$\s{ba}$}
\relabel{X2}{$\s{X}$}
\relabel{X1}{$\s{X}$}
\relabel{X3}{$\s{X}$}
\relabel{X4}{$\s{X}$}
\relabel{X5}{$\s{X}$}
\relabel{Y1}{$\s{Y}$}
\relabel{Y2}{$\s{Y}$}
\relabel{Y3}{$\s{Y}$}
\relabel{Y4}{$\s{Y}$}
\relabel{Y5}{$\s{Y}$}
\relabel{R41}{${R4}$}
\relabel{R42}{${R4}$}
\relabel{R43}{${R4}$}
\relabel{R21}{${R2}$}
\relabel{R22}{${R2}$}
\endrelabelbox}
\caption{A commutation relation.}
\label{com2}
\end{figure}

\begin{figure}
\centerline{\relabelbox 
\epsfysize 8cm
\epsfbox{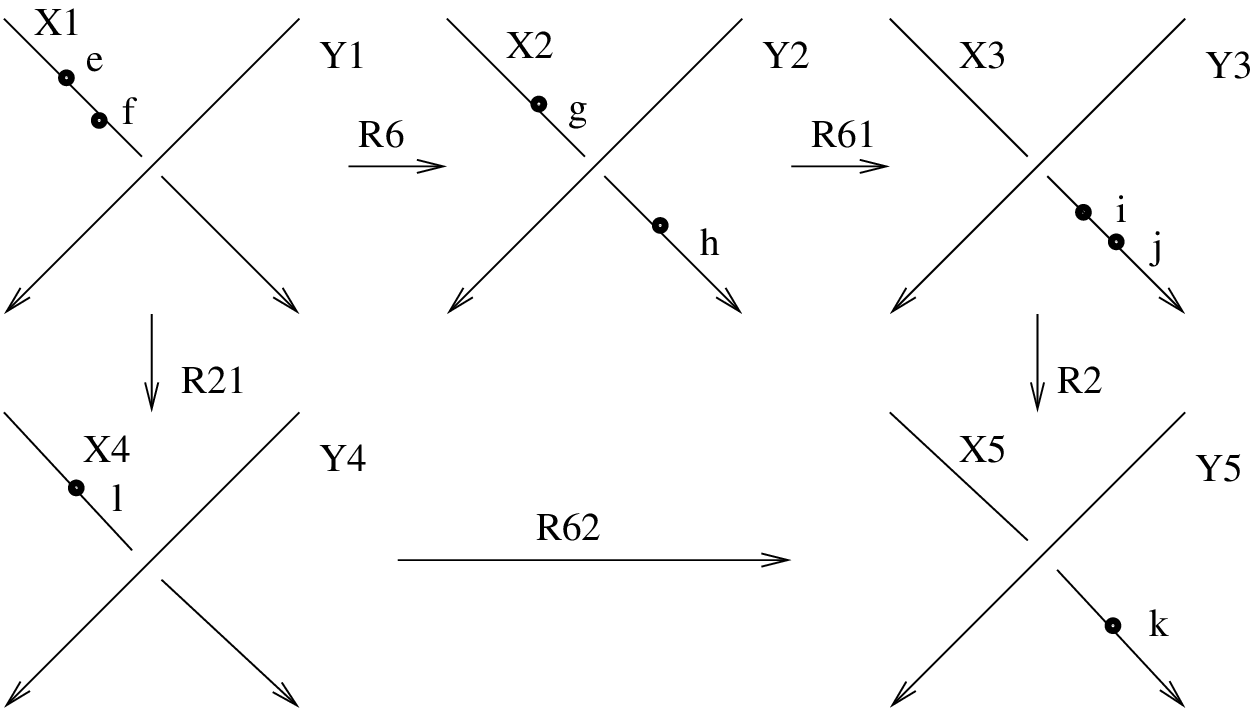}
\relabel{e}{$\s{a}$}
\relabel{f}{$\s{b}$}
\relabel{g}{$\s{a}$}
\relabel{h}{$\s{Y^{-1}\t b}$}
\relabel{i}{$\s{Y^{-1}\t a}$}
\relabel{j}{$\s{Y^{-1}\t b}$}
\relabel{k}{$\s{Y^{-1}(ba)}$}
\relabel{l}{$\s{ba}$}
\relabel{X1}{$\s{X}$}
\relabel{X2}{$\s{X}$}
\relabel{X3}{$\s{X}$}
\relabel{X4}{$\s{X}$}
\relabel{X5}{$\s{X}$}
\relabel{Y1}{$\s{Y}$}
\relabel{Y2}{$\s{Y}$}
\relabel{Y3}{$\s{Y}$}
\relabel{Y4}{$\s{Y}$}
\relabel{Y5}{$\s{Y}$}
\relabel{R6}{$R6$}
\relabel{R61}{$R6$}
\relabel{R62}{$R6$}
\relabel{R2}{$R2$}
\relabel{R21}{$R2$}
\endrelabelbox}
\caption{A commutation relation.}
\label{com3}
\end{figure}

The relations $R1$ to $R6$ commute with each other in the sense presented in
figures \ref{com1}, \ref{com2}, \ref{com3}. The calculations on the pictures have obvious counterparts for diagrams with the opposite crossing in the middle. 

\subsection{Reidemeister Moves}

We now prove that if $D$ and $D'$ are related by Reidemeister moves, then there are naturally associated linear maps $\V(D) \to \V(D')$. If $m$ is such a move, call this linear map  $F(m)$.

\subsubsection{Reidemeister-I}
\begin{figure}
\centerline{\relabelbox 
\epsfysize 3cm
\epsfbox{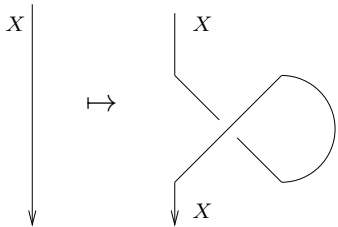} 
\relabel{t}{$\mapsto$}
\relabel{X1}{$\s{X}$}
\relabel{X2}{$\s{X}$}
\relabel{X3}{$\s{X}$}
\endrelabelbox }
\caption{Map associated to positive Reidemeister-I move.}
\label{MapR1}
\end{figure}

\begin{figure}
\centerline{\relabelbox 
\epsfysize 3cm
\epsfbox{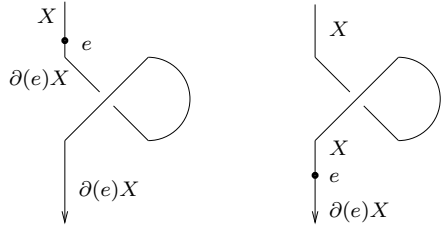} 
\relabel{X1}{$\s{X}$}
\relabel{X2}{$\s{\d(e)X}$}
\relabel{X3}{$\s{\d(e)X}$}
\relabel{X4}{$\s{X}$}
\relabel{X5}{$\s{X}$}
\relabel{X6}{$\s{\d(e)X}$}
\relabel{e}{$\s{e}$}
\relabel{f}{$\s{e}$}
\endrelabelbox }
\caption{Identity needing proof.}
\label{id1}
\end{figure}

\begin{figure}
\centerline{\relabelbox 
\epsfysize 3cm
\epsfbox{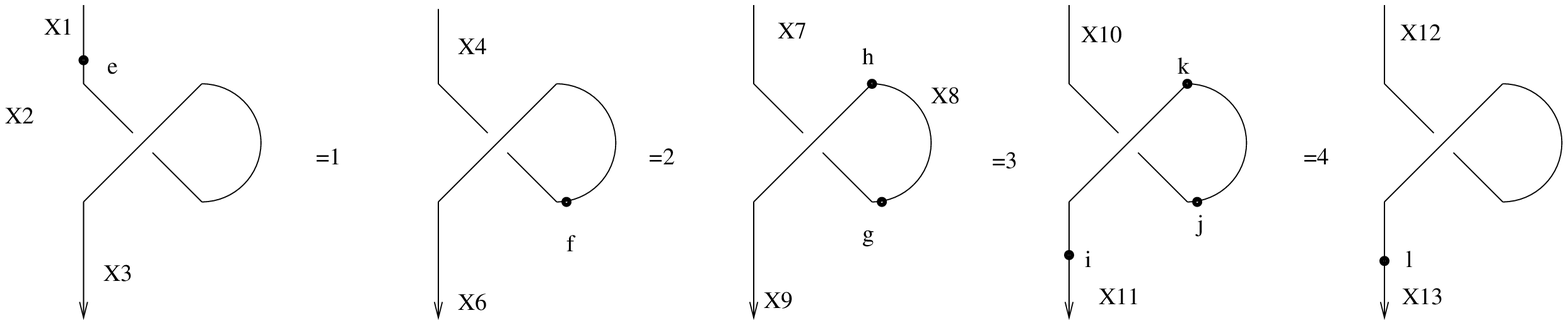} 
\relabel{X8}{$\s{X}$}
\relabel{X9}{$\s{\d(e)X}$}
\relabel{X10}{$\s{X}$}
\relabel{X11}{$\s{\d(e)X}$}
\relabel{X12}{$\s{X}$}
\relabel{X13}{$\s{\d(e)X}$}
\relabel{X1}{$\s{X}$}
\relabel{X2}{$\s{\d(e)X}$}
\relabel{X3}{$\s{\d(e)X}$}
\relabel{X4}{$\s{X}$}
\relabel{X6}{$\s{\d(e)X}$}
\relabel{X7}{$\s{X}$}
\relabel{e}{$\s{e}$}
\relabel{f}{$\s{X^{-1} \t e}$}
\relabel{g}{$\s{e^{-1}X^{-1} \t e }$}
\relabel{h}{$\s{e}$}
\relabel{i}{$\s{e}$}
\relabel{j}{$\s{X^{-1}\t e^{-1} e}$}
\relabel{k}{$\s{e^{-1}X^{-1} \t e}$}
\relabel{l}{$\s{e}$}
\relabel{=1}{$=$}
\relabel{=2}{$=$}
\relabel{=3}{$=$}
\relabel{=4}{$=$}
\endrelabelbox }
\caption{Proof of the identity in figure \ref{id1}.}
\label{id1proof}
\end{figure}

\begin{figure}
\centerline{\relabelbox 
\epsfysize 7cm
\epsfbox{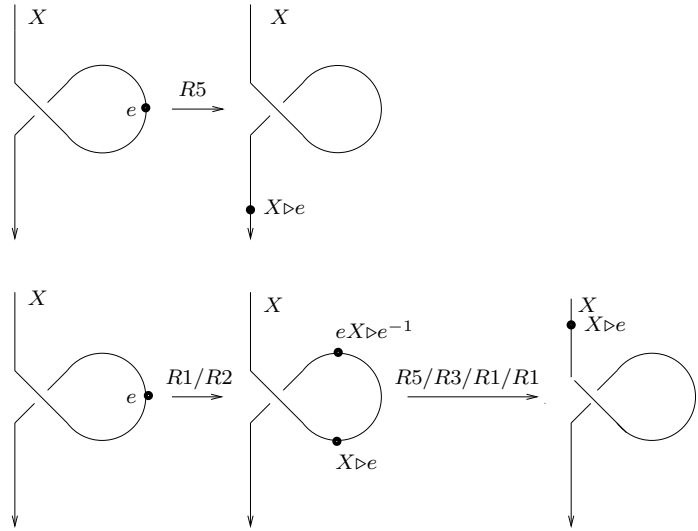} 
\relabel{X}{$\s{X}$}
\relabel{Y}{$\s{X}$}
\relabel{Z}{$\s{X}$}
\relabel{W}{$\s{X}$}
\relabel{T}{$\s{X}$}
\relabel{a}{$\s{e}$}
\relabel{b}{$\s{X \t e}$}
\relabel{c}{$\s{e}$}
\relabel{d}{$\s{e X \t e^{-1}}$}
\relabel{e}{$\s{X \t e}$}
\relabel{f}{$\s{X \t e}$}
\relabel{R1}{$\s{R5}$}
\relabel{R2}{$\s{R1/R2}$}
\relabel{R3}{$\s{R5/R3/R1/R1}$}
\endrelabelbox }
\caption{Two identities.}
\label{R1neg}
\end{figure}

As pointed out in the second chapter, since we are considering oriented knot diagrams, there
are four different cases of the  Reidemeister-I move, considered up to planar
isotopy. 
Suppose $D'$ is related to $D$ by a positive Reidemeister-I move $m$.
Positive means that it transforms a straight strand into a kink. We define a map $\V(D)\ra{F(m)} \V(D')$ as in figure \ref{MapR1}. Only one kind of Reidemeister-I is shown, but the other cases are perfectly analogous.
To prove $F(m)$ is well defined, we need to prove the equality of figure
\ref{id1}. This is done in figure \ref{id1proof}. The proof for the other
cases of Reidemeister-I move is analogous. If $m$ is a negative Reidemeister-I move, then the map $F(m)$
associated to it is also defined from figure \ref{MapR1}. To prove that the
definition of  $F(m)$ is correct, we need to use the identities of figure \ref{R1neg} together with
the one of
figure \ref{id1}. Note that we need to consider their analogues for the other types of Reidemeister-I move. The identities of figure \ref{R1neg} should be interpreted in
the light of \ref{commutation}.  It is important to note that if $m$ is a Reidemeister-I move we always have:  $F(m^{-1})=F(m)^{-1}$. 

\begin{Exercise}\label{Refer4}
Use Remark \ref{argument} to give a simple proof of all identities used to
prove that the map $F(m)$, where $m$ is a Reidemeister-I move,  is well defined.\end{Exercise}
\subsubsection{Reidemeister-II} 
\begin{figure}
\centerline{\relabelbox 
\epsfysize 3cm
\epsfbox{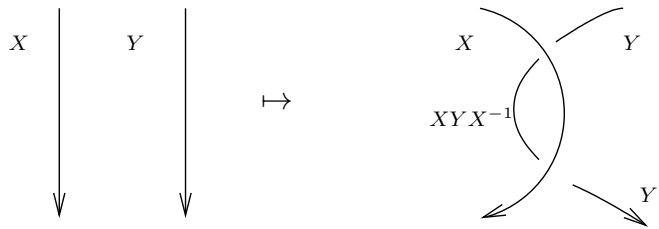} 
\relabel{X1}{$\s{X}$}
\relabel{X2}{$\s{X}$}
\relabel{Y1}{$\s{Y}$}
\relabel{Y2}{$\s{Y}$}
\relabel{Y3}{$\s{Y}$}
\relabel{Z}{$\s{XYX^{-1}}$}
\relabel{m}{$\mapsto$}
\endrelabelbox }
\caption{Map assigned to Reidemeister-II.}
\label{MapR2}
\end{figure}

\begin{figure}
\centerline{\relabelbox 
\epsfysize 3cm
\epsfbox{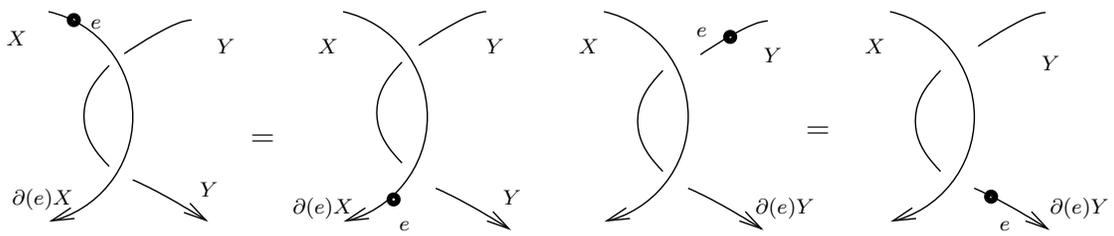} 
\relabel{X1}{$\s{X}$}
\relabel{X2}{$\s{X}$}
\relabel{X3}{$\s{\d(e)X}$}
\relabel{X4}{$\s{\d(e)X}$}
\relabel{X5}{$\s{X}$}
\relabel{X6}{$\s{X}$}
\relabel{Y1}{$\s{Y}$}
\relabel{Y8}{$\s{Y}$}
\relabel{Y2}{$\s{Y}$}
\relabel{Y3}{$\s{Y}$}
\relabel{Y4}{$\s{Y}$}
\relabel{Y5}{$\s{\d(e)Y}$}
\relabel{Y6}{$\s{Y}$}
\relabel{Y7}{$\s{\d(e)Y}$}
\relabel{=1}{$=$}
\relabel{=2}{$=$}
\relabel{e}{$\s{e}$}
\relabel{f}{$\s{e}$}
\relabel{g}{$\s{e}$}
\relabel{h}{$\s{e}$}
\endrelabelbox }
\caption{Two identities that need proving.}
\label{id2}
\end{figure}

\begin{figure}
\centerline{\relabelbox 
\epsfysize 4cm
\epsfbox{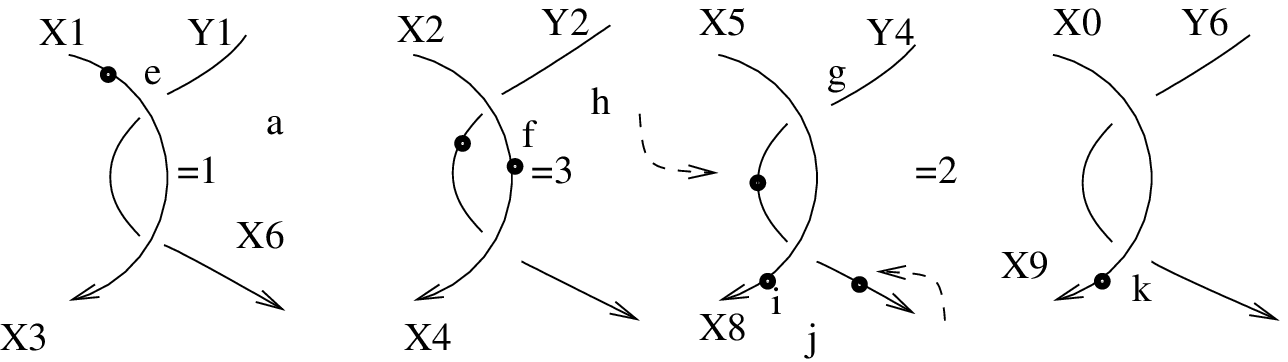} 
\relabel{X0}{$\s{X}$}
\relabel{X1}{$\s{X}$}
\relabel{X3}{$\s{\d(e)X}$}
\relabel{Y1}{$\s{Y}$}
\relabel{Y2}{$\s{Y}$}
\relabel{X2}{$\s{X}$}
\relabel{X4}{$\s{\d(e)X}$}
\relabel{X8}{$\s{\d(e)X}$}
\relabel{a}{$\s{e XYX^{-1} \t e^{-1}}$}
\relabel{X9}{$\s{\d(e)X}$}
\relabel{e}{$\s{e}$}
\relabel{f}{$\s{e}$}
\relabel{h}{$\s{e XYX^{-1} \t e^{-1}}$}
\relabel{X6}{$\s{\d(e)XYX^{-1}\d(e^{-1})}$}
\relabel{X5}{$\s{X}$}
\relabel{Y4}{$\s{Y}$}
\relabel{Y6}{$\s{Y}$}
\relabel{i}{$\s{e}$}
\relabel{k}{$\s{e}$}
\relabel{j}{$\s{X^{-1}\t e^{-1}X^{-1} \d(e)XYX^{-1}\d(e^{-1})\t e}$}
\relabel{=1}{$=$}
\relabel{=2}{$=$}
\relabel{=3}{$=$}
\endrelabelbox}
\caption{Proof of identity in figure \ref{id2}.}
\label{id2proof}
\end{figure}

Let $D$ and $D'$ be knot diagrams differing by  Reidemeister move number II. Call it
 $m$. Define a map $\V(D) \ra{F(m)} \V(D')$ as in figure \ref{MapR2}. The
 definition of $F(m)$ for the other types of Reidemeister-II moves is
 analogous. To prove that $F(m)$ is well defined, we need to prove the
 equalities of figure \ref{id2}. The most difficult one is the first and the
 proof of it 
 appears in figure \ref{id2proof}.  The last identity follows from
\begin{align*}
X^{-1}\t e^{-1} &X^{-1}\d(e)XYX^{-1}\d(e^{-1})\t e X^{-1}\t(eXYX^{-1}\t e^{-1}
)\\
&=X^{-1}\t e^{-1} X^{-1}\t e YX^{-1}\t e X^{-1} \t e^{-1} X^{-1}\t e YX^{-1}\t e^{-1}
\\
&=X^{-1}\t e^{-1} X^{-1}\t e YX^{-1}\t e YX^{-1}\t e^{-1}
\\
&=1_E
\end{align*}

\begin{figure}
\centerline{\relabelbox 
\epsfysize 8cm
\epsfbox{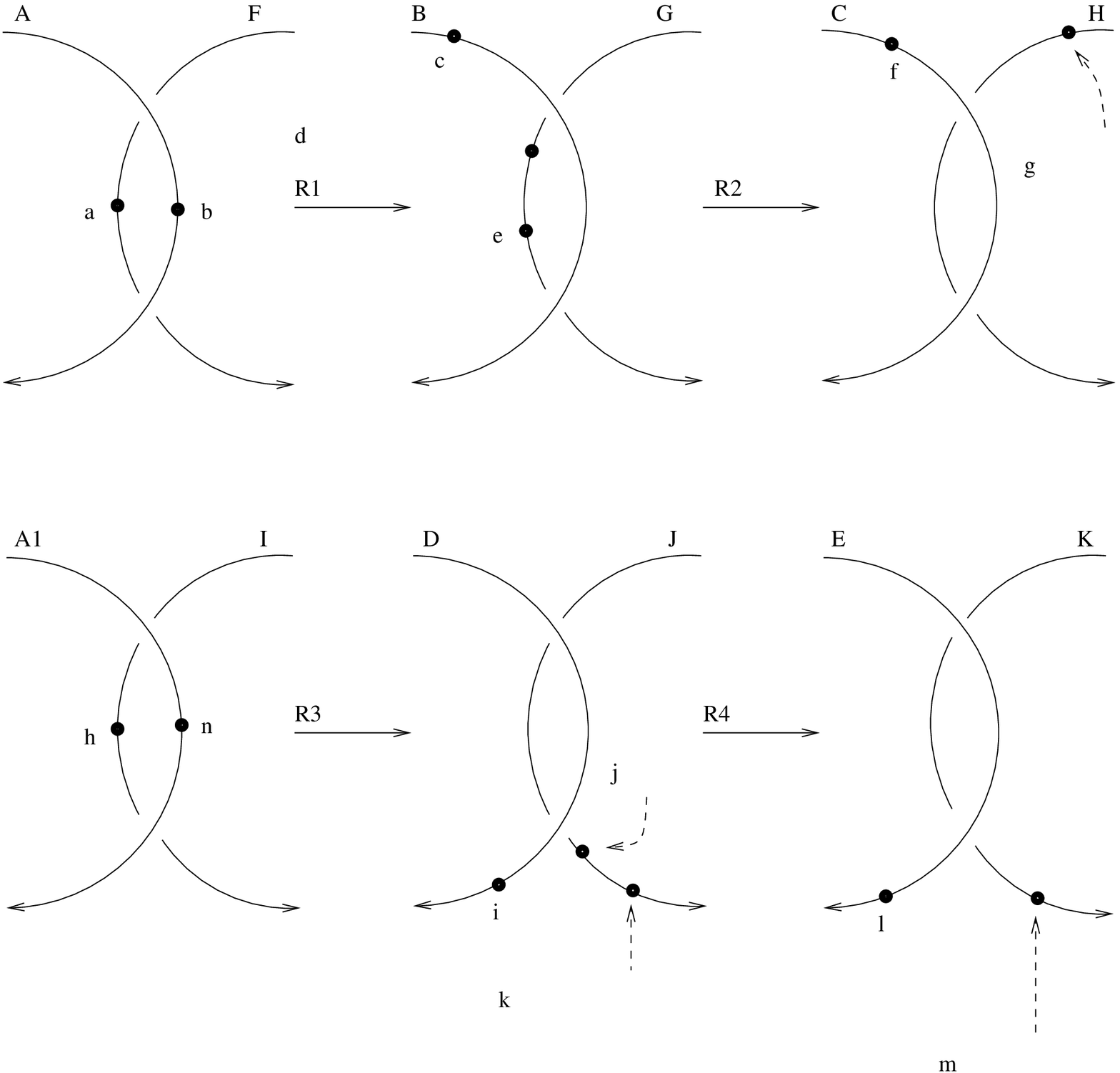}
\relabel{A}{$\s{X}$}
\relabel{B}{$\s{X}$}
\relabel{C}{$\s{X}$}
\relabel{D}{$\s{X}$}
\relabel{E}{$\s{X}$}
\relabel{A1}{$\s{X}$}
\relabel{F}{$\s{Y}$}
\relabel{G}{$\s{Y}$}
\relabel{H}{$\s{Y}$}
\relabel{I}{$\s{Y}$}
\relabel{J}{$\s{Y}$}
\relabel{K}{$\s{Y}$}
\relabel{a}{$\s{a}$}
\relabel{b}{$\s{b}$}
\relabel{c}{$\s{b}$}
\relabel{d}{$\s{XYX^{-1}\t b b ^{-1}}$}
\relabel{e}{$\s{a}$}
\relabel{f}{$\s{b}$}
\relabel{g}{$\s{X^{-1} \t b^{-1} X^{-1}\t a YX^{-1}\t b }$}
\relabel{h}{$\s{a}$}
\relabel{n}{$\s{b}$}
\relabel{i}{$\s{b}$}
\relabel{j}{$\s{X^{-1} \t b^{-1} YX^{-1}\t b }$}
\relabel{k}{$\s{X^{-1}\t b^{-1} X^{-1}\t a X^{-1}\t b}$}
\relabel{l}{$\s{b}$}
\relabel{m}{$\s{X^{-1}\t b^{-1} X^{-1}\t a YX^{-1} \t b }$}
\relabel{R1}{$\s{R3}$}
\relabel{R2}{$\s{R2/R5}$}
\relabel{R3}{$\s{R4/R6}$}
\relabel{R4}{$\s{R2}$}
\endrelabelbox}
\caption{Two identities.}
\label{R2neg}
\end{figure}

If $m$ is a negative Reidemeister-II move, then the map $F(m)$ is similarly defined
from figure \ref{MapR2}. To prove that the definition of $F(m)$ makes sense, we need to use the
identities of figure \ref{R2neg}. These are not ambiguous due to
\ref{commutation}. Note that we have $F(m^{-1})=F(m)^{-1}$ if $m$ is
 a Reidemeister-II move. The other types of Reidemeister-II move are dealt with similarly. The identities which we used  to prove that the
maps $F(m):\V(D) \to \V(D')$, for $m$ a Reidemeister-II move,  are well defined can also be shown using remark
\ref{argument} (c.f. exercise \ref{Refer4}). 
\subsubsection{Reidemeister-III}
\begin{figure}
\centerline{\relabelbox 
\epsfysize 4cm
\epsfbox{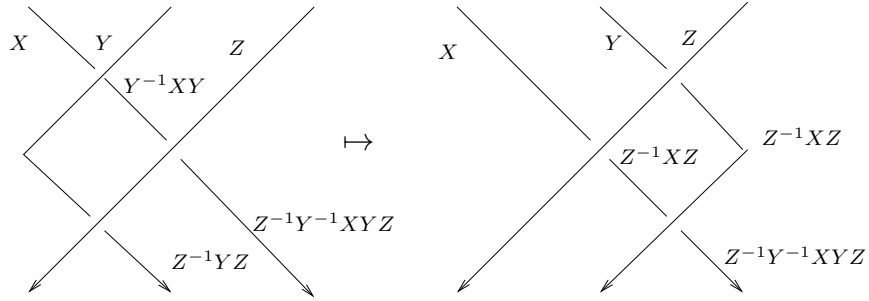}
\relabel{t}{$\mapsto$}
 \relabel{X}{$\s{X}$}
\relabel{Y}{$\s{Y}$}
\relabel{Z}{$\s{Z}$}
\relabel{X1}{$\s{Y^{-1}XY}$}
\relabel{Y1}{$\s{Z^{-1}YZ}$}
\relabel{X2}{$\s{Z^{-1}Y^{-1}XYZ}$}
\relabel{A}{$\s{X}$}
\relabel{B}{$\s{Y}$}
\relabel{C}{$\s{Z}$}
\relabel{B1}{$\s{Z^{-1}XZ}$}
\relabel{A2}{$\s{Z^{-1}XZ}$}
\relabel{A3}{$\s{Z^{-1}Y^{-1}XYZ}$}
\endrelabelbox}
\caption{Map assigned to Reidemeister-III.}
\label{MapR3}
\end{figure}

\begin{figure}
\centerline{\relabelbox 
\epsfysize 10cm
\epsfbox{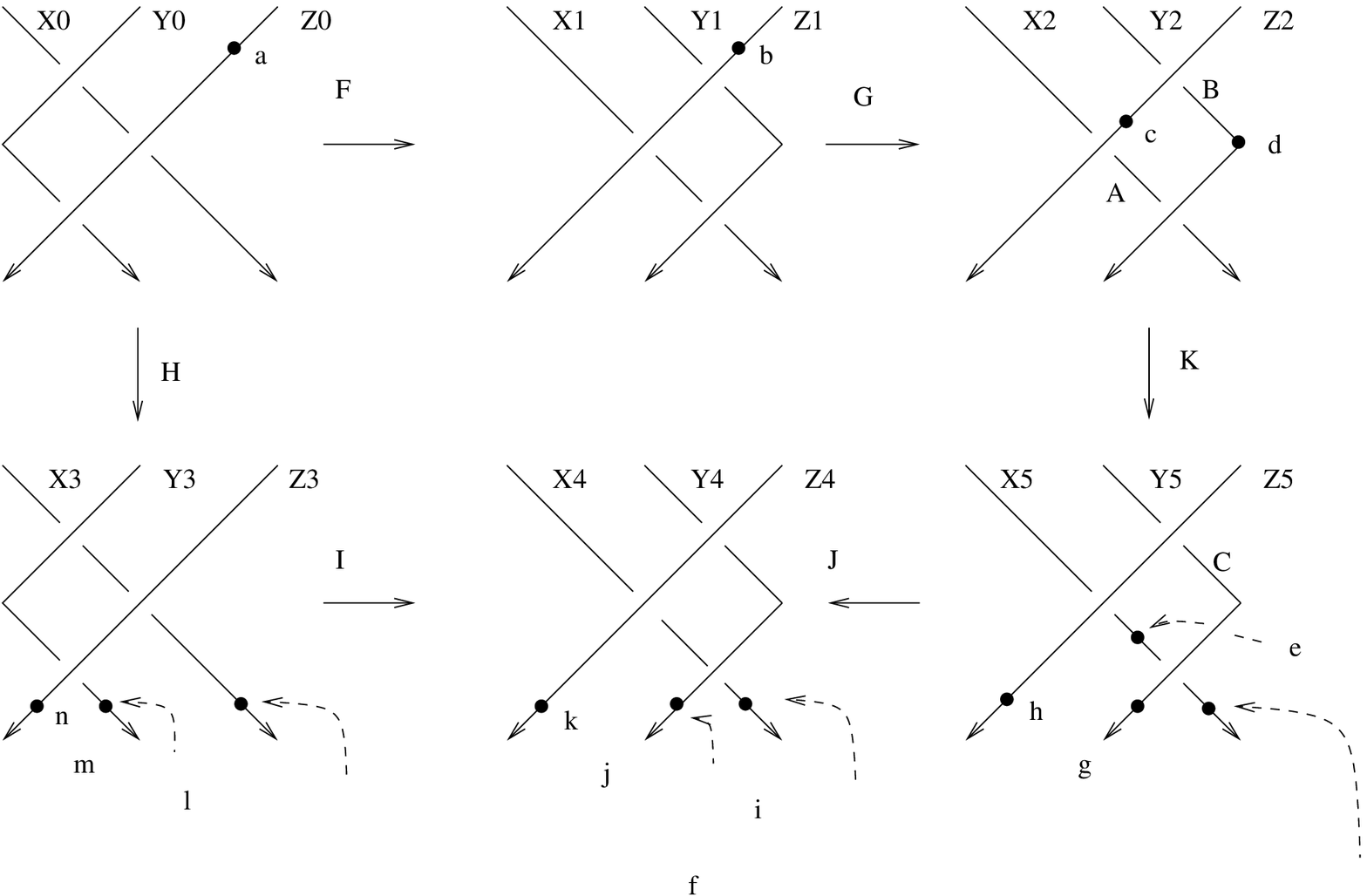}
\relabel{X0}{$\s{X}$}
\relabel{X1}{$\s{X}$}
\relabel{X2}{$\s{X}$}
\relabel{X3}{$\s{X}$}
\relabel{X4}{$\s{X}$}
\relabel{X5}{$\s{X}$}
\relabel{Y0}{$\s{Y}$}
\relabel{Y1}{$\s{Y}$}
\relabel{Y2}{$\s{Y}$}
\relabel{Y3}{$\s{Y}$}
\relabel{Y4}{$\s{Y}$}
\relabel{Y5}{$\s{Y}$}
\relabel{Z0}{$\s{Z}$}
\relabel{Z1}{$\s{Z}$}
\relabel{Z2}{$\s{Z}$}
\relabel{Z3}{$\s{Z}$}
\relabel{Z4}{$\s{Z}$}
\relabel{Z5}{$\s{Z}$}
\relabel{A}{$\s{A}$}
\relabel{B}{$\s{Z^{-1}YZ}$}
\relabel{C}{$\s{Z^{-1}YZ}$}
\relabel{a}{$\s{e}$} 
\relabel{b}{$\s{e}$}
\relabel{c}{$\s{e}$}
\relabel{d}{$\s{Z^{-1}\t e^{-1}Z^{-1}Y \t e } $}
\relabel{e}{$\s{Z^{-1}\t  e^{-1} Z^{-1}X \t e }$}
\relabel{f}{$\s{Z^{-1}\t e^{-1} Z^{-1}Y^{-1}\t e  Z^{-1}Y^{-1}\d(e^{-1})X \d(e) \t e^{-1} Z^{-1}Y^{-1}\d(e^{-1})X\d(e)Y \t e}$}
\relabel{g}{$\s{Z^{-1} \t e^{-1} Z^{-1}Y \t e}$}
\relabel{h}{$\s{e}$}
\relabel{i}{$\s{Z^{-1} \t e^{-1}Z^{-1} Y^{-1} XY \t e}$}
\relabel{j}{$\s{Z^{-1}\t e^{-1} Z^{-1}Y \t e}$}
\relabel{k}{$\s{e}$}
\relabel{l}{$\s{Z^{-1} \t e^{-1}Z^{-1} Y^{-1} XY \t e}$}
\relabel{m}{$\s{Z^{-1}\t e^{-1} Z^{-1}Y \t e}$}
\relabel{n}{$\s{e}$}
\relabel{F}{$F(m)$}
\relabel{G}{$R4$}
\relabel{K}{$R4$}
\relabel{J}{$R6/R2$}
\relabel{H}{$R4$}
\relabel{I}{$F(m)$}
\endrelabelbox}
\caption{Proof that $F(m)$ is well defined if $m$ is a Reidemeister-III move. Here $A=Z^{-1}\d(e^{-1})X \d(e)Z$.}
\label{PROOF}
\end{figure}

As we mentioned in the second chapter, since we are considering oriented knot diagrams, there are sixteen versions of
the Reidemeister move number III. The positive move has the direction indicated in the
figure \ref{R3}. Let
$D$ and $D'$ differ by a Reidemeister-III move $m$. Suppose we have a dotting
of $D$. We can always move the vertices away from the area in question. For
definiteness, suppose we move all vertices in the direction defined by the
orientation of the edges, starting with the top strand and finishing with the 
bottom one. This operation verifies the relations $R1$ to $R6$ due to the
commutation relation presented in \ref{commutation}. Having done these
changes, the definition of $\V(D)\ra{F(m)} \V(D')$ is straightforward and appears in
figure \ref{MapR3}. It is analogous for the other types of Reidemeister-III
move. To prove $F(m)$ is well defined we still need to prove a set of relations similar to the relations in figures \ref{id1} and \ref{id2}, and then apply the results of \ref{commutation}. The proof of the most difficult one appears in figure \ref{PROOF}. The last equality follows from:
\begin{multline*}
(Z^{-1}\t e^{-1} Z^{-1}Y^{-1}\t e  Z^{-1}Y^{-1}\d(e^{-1})X\d(e)\t e^{-1} Z^{-1}Y^{-1}\d(e^{-1})X\d(e)Y \t e)\\ Z^{-1}Y^{-1}Z \t (Z^{-1}\t  e^{-1} Z^{-1}X \t e)\\
=Z^{-1}\t e^{-1}
 Z^{-1}Y^{-1}\t e
 Z^{-1}Y^{-1} \t e^{-1} 
Z^{-1}Y^{-1}X \t e^{-1} Z^{-1} Y ^{-1} \t e \\ Z^{-1}Y^{-1} \t e^{-1} Z^{-1} Y^{-1} X \t e Z^{-1}Y^{-1}  X Y \t e Z^{-1}Y^{-1}  X\t e^{-1} Z^{-1} Y^{-1} \t e  \\
Z^{-1}Y^{-1} \t  e^{-1}   Z^{-1}Y^{-1} X \t e\\
=Z^{-1}\t e^{-1}   Z^{-1}Y^{-1}  X Y \t e.
\end{multline*}
To prove that the  maps associated to the other cases of Reidemeister-III move
are well defined we can proceed similarly. The definition of the map $F(m)$ if
$m$ is a negative Reidemeister-III move is totally analogous. We always have $F(m^{-1})=F(m)^{-1}$. It is possible to
prove that all maps associated to Reidemeister-III move are well defined from  remark \ref{argument} and \ref{commutation}, as it was for the
other Reidemeister moves. 

\subsection{Morse Moves}
Similarly to Reidemeister moves, if $D$ and $D'$ are knot diagrams related by a  Morse move $m$, then there is  naturally associated a map  $\V(D) \ra{F(m)} \V(D')$. 
\subsubsection{Saddle Points}
\begin{figure}
\centerline{\relabelbox 
\epsfysize 2cm
\epsfbox{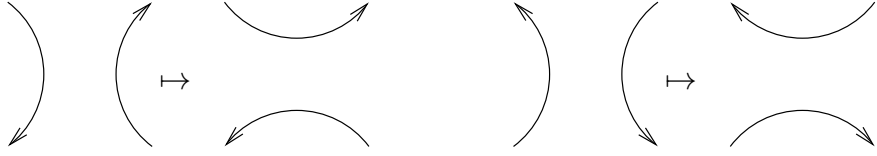} 
\relabel{t1}{$\mapsto$}
\relabel{t2}{$\mapsto$}
\endrelabelbox}
\caption{Oriented saddle point move.}
\label{saddle}
\end{figure}

\begin{figure}
\centerline{\relabelbox 
\epsfysize 2cm
\epsfbox{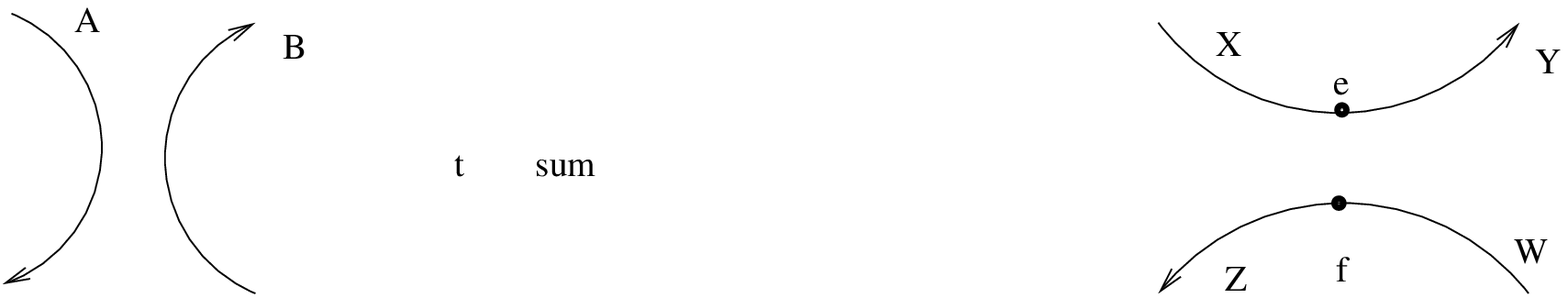}
\relabel{t}{$\mapsto$}
\relabel{A}{$\s{X}$}
\relabel{B}{$\s{Y}$}
\relabel{X}{$\s{X}$}
\relabel{Y}{$\s{\d(e)X}$}
\relabel{Z}{$\s{X}$}
\relabel{W}{$\s{\d(e)X}$}
\relabel{e}{$\s{e}$}
\relabel{f}{$\s{e^{-1}}$}
\relabel{sum}{$\displaystyle{\frac{1}{\# E}{\sum_{e \in E} \delta(Y,\d(e)X)}}$}
\endrelabelbox}
\caption{Map associated to saddle point moves.}
\label{Mapsaddle}
\end{figure}

\begin{figure}
\centerline{\relabelbox 
\epsfysize 2.5cm
\epsfbox{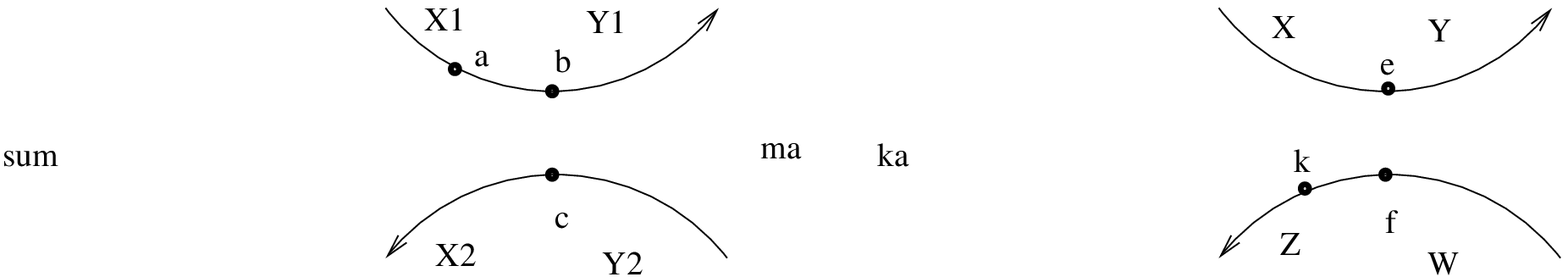}
\relabel{ma}{$=$}
\relabel{ka}{$\displaystyle{{\sum_{e \in E} \delta(Y,\d(e)X)}}$}
\relabel{X1}{$\s{X}$}
\relabel{Y1}{$\s{\d(ef) X}$}
\relabel{X2}{$\s{\d(f)X}$}
\relabel{Y2}{$\s{\d(ef)X}$}
\relabel{X}{$\s{X}$}
\relabel{Y}{$\s{\d(e)X}$}
\relabel{Z}{$\s{\d(f)X}$}
\relabel{Y}{$\s{\d(e)X}$}
\relabel{W}{$\s{\d(e)X}$}
\relabel{a}{$\s{f}$}
\relabel{b}{$\s{e}$}
\relabel{c}{$\s{e^{-1}}$}
\relabel{e}{$\s{e}$}
\relabel{f}{$\s{e^{-1}}$}
\relabel{k}{$\s{f}$}
\relabel{sum}{$\displaystyle{{\sum_{e \in E} \delta(Y,\d(e)\d(f)X)}}$}
\endrelabelbox}
\caption{An identity.}
\label{id4}
\end{figure}

\begin{figure}
\centerline{\relabelbox
\epsfysize 3cm
\epsfbox{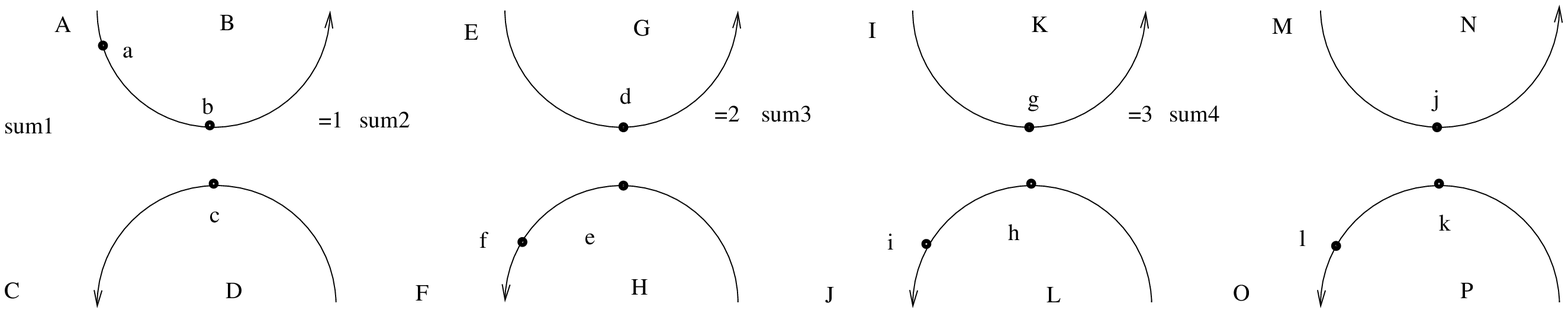}
\relabel{A}{$\s{X}$}
\relabel{E}{$\s{X}$}
\relabel{I}{$\s{X}$}
\relabel{M}{$\s{X}$}
\relabel{B}{$\s{\d(fe)X}$}
\relabel{G}{$\s{\d(fe)X}$}
\relabel{K}{$\s{\d(fe)X}$}
\relabel{D}{$\s{\d(fe)X}$}
\relabel{H}{$\s{\d(fe)X}$}
\relabel{L}{$\s{\d(fe)X}$}
\relabel{C}{$\s{\d(e)X}$}
\relabel{F}{$\s{\d(e)X}$}
\relabel{J}{$\s{\d(e)X}$}
\relabel{O}{$\s{\d(e)X}$}
\relabel{P}{$\s{\d(g)X}$}
\relabel{N}{$\s{\d(g)X}$}
\relabel{a}{$\s{f}$}
\relabel{b}{$\s{e}$}
\relabel{c}{$\s{e^{-1}}$}
\relabel{d}{$\s{ef}$}
\relabel{e}{$\s{f^{-1}e^{-1}}$}
\relabel{f}{$\s{f}$}
\relabel{g}{$\s{ef}$}
\relabel{h}{$\s{(ef)^{-1}}$}
\relabel{i}{$\s{f}$}
\relabel{j}{$\s{g}$}
\relabel{k}{$\s{g^{-1}}$}
\relabel{l}{$\s{f}$}
\relabel{=1}{$=$}
\relabel{=2}{$=$}
\relabel{=3}{$=$}
\relabel{sum1}{$\displaystyle{\sum_{\substack{e \in E\\ \d(fe)X=Y}}}$}
\relabel{sum2}{$\displaystyle{\sum_{\substack{e \in E\\ \d(fe)X=Y}}}$}
\relabel{sum3}{$\displaystyle{\sum_{\substack{e \in E\\ \d(fe)X=Y}}}$}
\relabel{sum4}{$\displaystyle{\sum_{\substack{g \in E\\ \d(g)X=Y}}}$}
\endrelabelbox}
\caption{Proof of the identity in figure \ref{id4}.}
\label{id4proof}
\end{figure}

There are two kinds of saddle point moves in the oriented case. They are
depicted in figure \ref{saddle}.
Let $D$ and $D'$ differ by a  saddle point move. We can define a map $\V(D) \to
\V(D')$. We make the definition for the first move, since it is analogous for
the second one. It appears in figure \ref{Mapsaddle}. The map is well defined due to
the identity in figure \ref{id4}. This identity is proved in figure \ref{id4proof}.

\subsubsection{Cups and Caps}
\begin{figure}
\centerline{\relabelbox
\epsfysize 1.5cm
\epsfbox{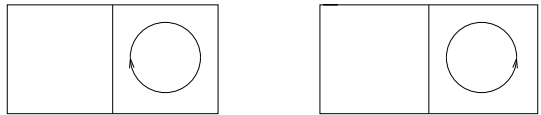}\endrelabelbox}
\caption{Births of a circle.}
\label{cup}
\end{figure}

\begin{figure}
\centerline{\relabelbox
  \epsfysize 1.5cm
\epsfbox{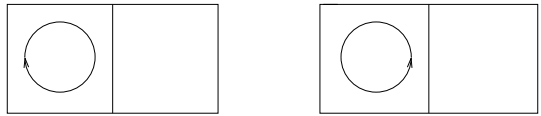}\endrelabelbox}
\caption{Deaths of a circle.}
\label{cap}
\end{figure}

\begin{figure}
\centerline{\relabelbox 
\epsfysize 1.5cm
\epsfbox{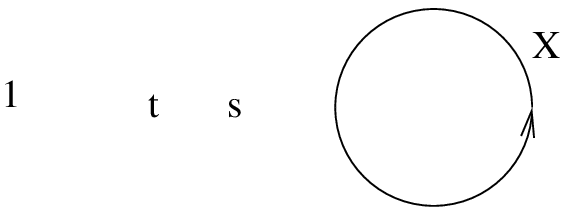}
\relabel{1}{$1$}
\relabel{t}{$\mapsto$}
\relabel{s}{$\displaystyle{\sum_{X \in G}}$}
\relabel{X}{$\s{X}$}
\endrelabelbox}
\caption{Map associated with births of a circle.}
\label{Mapcap}
\end{figure}

\begin{figure}
\centerline{\relabelbox 
\epsfysize 2.5cm
\epsfbox{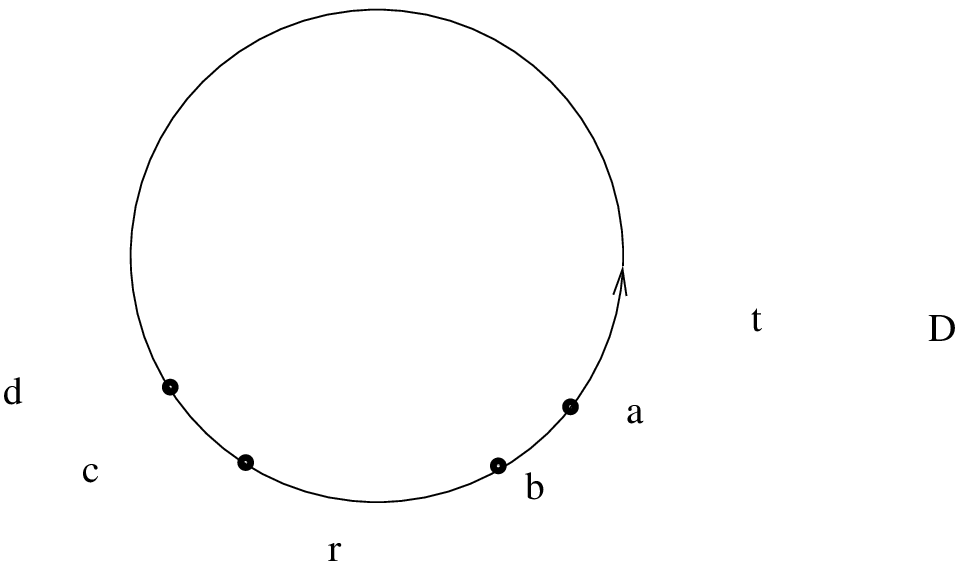}
\relabel{t}{$\mapsto$}
\relabel{a}{$\s{x_1}$}
\relabel{b}{$\s{x_2}$}
\relabel{c}{$\s{x_{n-1}}$}
\relabel{d}{$\s{x_n}$}
\relabel{r}{$\dots$}
\relabel{D}{$\# E \delta(x_1x_2...x_{n-1}x_n,1_E)$}
\endrelabelbox}
\caption{Map associated with deaths of a circle.}
\label{Mapcup}
\end{figure}

Since we are considering oriented diagrams, there are two kinds of births/deaths
of a circle. They are described in figures \ref{cup} and \ref{cap}. 
The maps assigned to births and deaths of circles  are described in figures  \ref{Mapcup} and \ref{Mapcap}. There, $1$ simply means what is left of the diagram. The other orientations are similar. It is easy to conclude that both types of maps are well defined. 

\subsection{Invariance}
If we are given a movie of a knotted surface $\S$, we can use the maps defined in
the previous section to give an element $I^4_\G(\S)$ of the ground field $\Q$. We assign the obvious map $\V(D) \to \V(D')$ if the knot diagram $D'$ is related to $D$ by a planar isotopy. To prove that it
is an invariant of knotted surfaces, we need to prove invariance under the
Movie Moves of Carter and Saito and interchanging distant critical points. The
invariance under interchanging distant critical points is trivial to verify
since all the maps defined are of a local nature. In the remainder of this
subsection, we prove invariance under the Movie Moves. We use the numbering of
the movie moves that appears in \cite{BN}. Note that, since we consider
oriented surfaces, we will need to consider not only the mirror image of the
all the movie moves but also all the possible orientations of the strands.
\subsubsection{Movie Moves $1$ to $5$}
The invariance under these movie moves is a consequence of the definition of
$F(m)$ for a Reidemeister move, since we always have $F(m^{-1})=F(m)^{-1}$. Note
that we can always suppose that no vertices get in the way, by moving them
away from the area in question and applying the commutation relations of \ref{commutation}.
\subsubsection{Movie Moves $6$ to $10$}\label{trivial}
We can always suppose that no vertices get in the way. The proof of invariance
is a straightforward calculation. The invariance under these movie moves is
known already from the fact that, not considering saddle points, we are simply
calculating the number of morphisms of the complement of a regular neighbourhood of the knotted surface $\S$ into $G$.
\subsubsection{Going Right Movie Move $11$}
\begin{figure}
\centerline{\relabelbox 
\epsfysize 1.5cm
\epsfbox{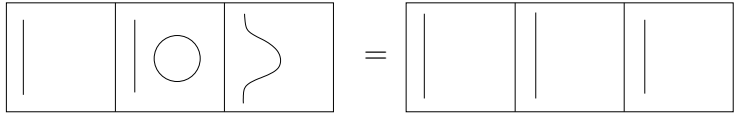}
\relabel{=}{$=$}
\endrelabelbox}
\caption{Movie Move 11.}
\label{MM11}
\end{figure}

\begin{figure}
\centerline{\relabelbox 
\epsfysize 2cm
\epsfbox{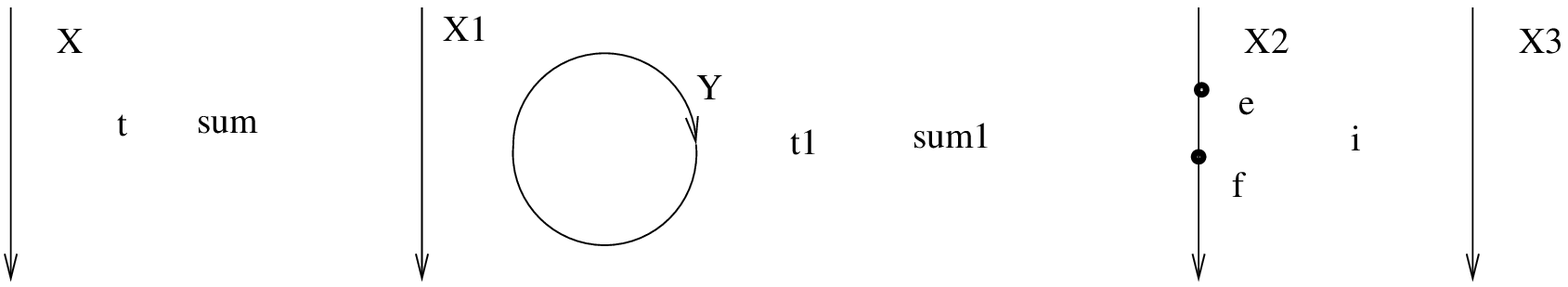}
\relabel{i}{$=$}
\relabel{t}{$\mapsto$}
\relabel{t1}{$\mapsto$}
\relabel{X}{$\s{X}$}
\relabel{X1}{$\s{X}$}
\relabel{X2}{$\s{X}$}
\relabel{X3}{$\s{X}$}
\relabel{Y}{$\s{Y}$}
\relabel{e}{$\s{e}$}
\relabel{f}{$\s{e^{-1}}$}
\relabel{sum}{$\displaystyle{\sum_{Y \in G}}$}
\relabel{sum1}{$\displaystyle{\frac{1}{\# E}\sum_{e \in E}}$}
\endrelabelbox}
\caption{Invariance under Going Right Movie Move $11$.}
\label{rightMM11}
\end{figure}

The movie move $11$ appears in figure \ref{MM11}. We need to consider all the
possible orientations. This movie move is reversible, so we consider the going
right and going left cases. The proof of invariance under the going right case  appears in figure
\ref{rightMM11}. We only consider one orientation, but the other case is
perfectly analogous.

\subsubsection{Going Left Movie Move $11$.}
\begin{figure}
\centerline{\relabelbox 
\epsfysize 2.3cm
\epsfbox{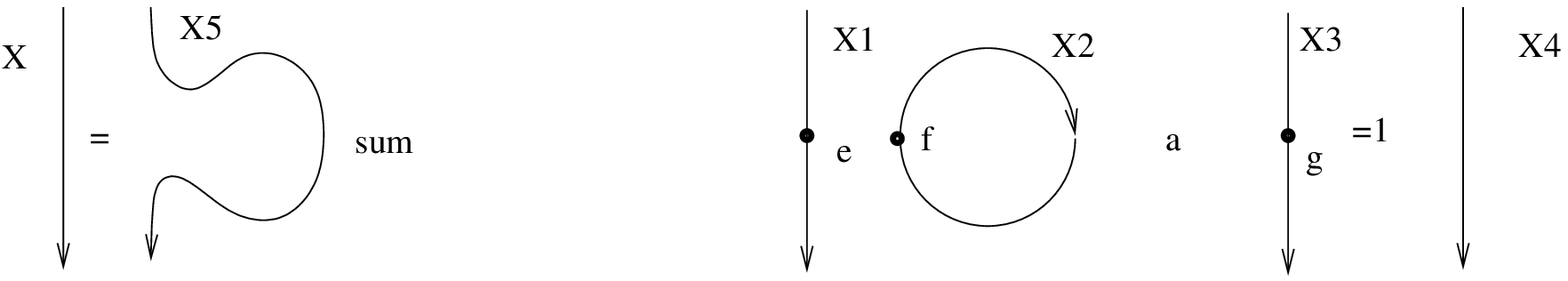}
\relabel{=}{$=$}
\relabel{=1}{$=$}
\relabel{X}{$\s{X}$}
\relabel{X1}{$\s{X}$}
\relabel{X2}{$\s{X}$}
\relabel{X3}{$\s{X}$}
\relabel{X4}{$\s{X}$}
\relabel{X5}{$\s{X}$}
\relabel{e}{$\s{e}$}
\relabel{f}{$\s{e^{-1}}$}
\relabel{g}{$\s{1_E}$}
\relabel{sum}{$\displaystyle{\mapsto \frac{1}{\# E}\sum_{e \in \mathrm{Ker}(\d)}}$}
\relabel{a}{$\mapsto$}
\endrelabelbox}
\caption{Invariance under Going Left Movie Move $11$.}
\label{leftMM11}
\end{figure}

The proof of invariance appears in figure \ref{leftMM11}.

\subsubsection{Movie Move 12}
\begin{figure}
\centerline{\relabelbox 
\epsfysize 1.5cm
\epsfbox{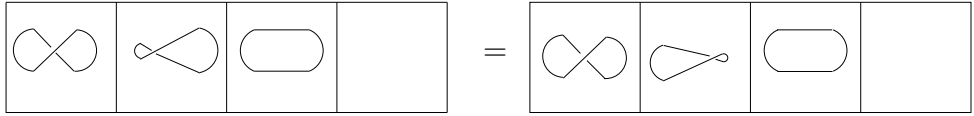}
\relabel{=}{$=$}
\endrelabelbox}
\caption{Movie Move 12.}
\label{MM12}
\end{figure}

\begin{figure}
\centerline{\relabelbox 
\epsfysize 3cm
\epsfbox{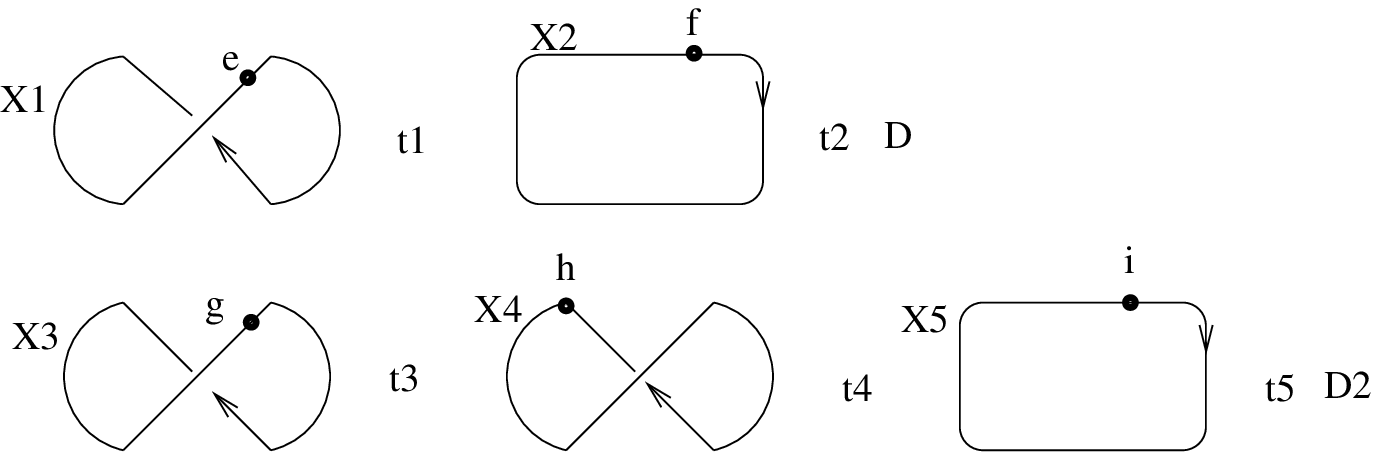}
\relabel{X1}{$\s{X}$}
\relabel{X2}{$\s{X}$}
\relabel{X3}{$\s{X}$}
\relabel{X4}{$\s{X}$}
\relabel{X5}{$\s{X}$}
\relabel{e}{$\s{e}$}
\relabel{f}{$\s{e}$}
\relabel{g}{$\s{e}$}
\relabel{h}{$\s{X ^{-1}\t e}$}
\relabel{i}{$\s{X^{-1} \t e}$}
\relabel{t1}{$\mapsto$}
\relabel{t2}{$\mapsto$}
\relabel{t3}{$\mapsto$}
\relabel{t4}{$\mapsto$}
\relabel{t5}{$\mapsto$}
\relabel{D}{$\#E \delta(e,1_E)$}
\relabel{D2}{$\#E \delta(X^{-1} \t e,1_E)$}
\endrelabelbox}
\caption{Proof of invariance under Going Right Movie Move $12$. }
\label{MM12proof}
\end{figure}

The movie move 12 appears in figure \ref{MM12}. We need to consider its mirror
image as well as a change on the  orientation. The going right move is trivial to
verify since there are no vertices involved. The invariance under going right
movie move 12 is a bit more complicated since some vertices can get in the
way. However, we can obviously suppose we have one only vertex, from relations $R1$ to $R6$; and the proof
of invariance for this case appears in figure \ref{MM12proof}. The arguments used to prove invariance under the mirror
images and the different orientations of this movie move are perfectly analogous.

\subsubsection{Movie Move 13}
\begin{figure}
\centerline{\relabelbox 
\epsfysize 2.5cm
\epsfbox{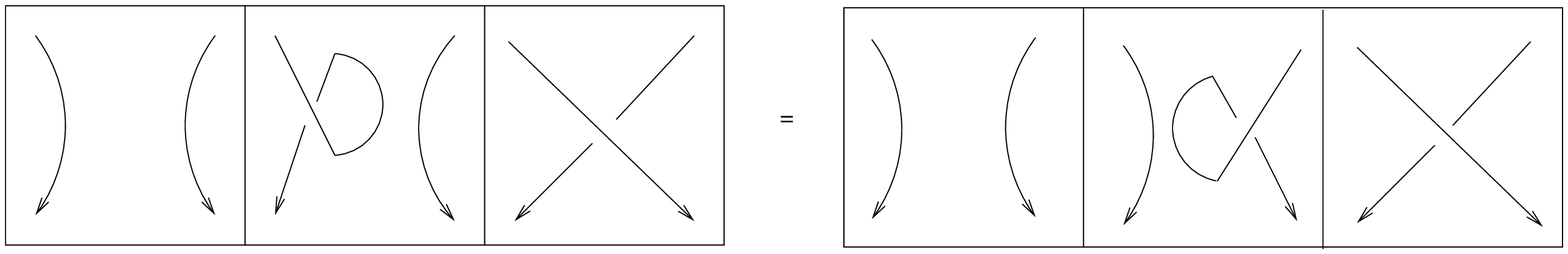}
\relabel{=}{$=$}
\endrelabelbox}
\caption{Movie Move 13.}
\label{MM13}
\end{figure}

\begin{figure}
\centerline{\relabelbox 
\epsfysize 4cm
\epsfbox{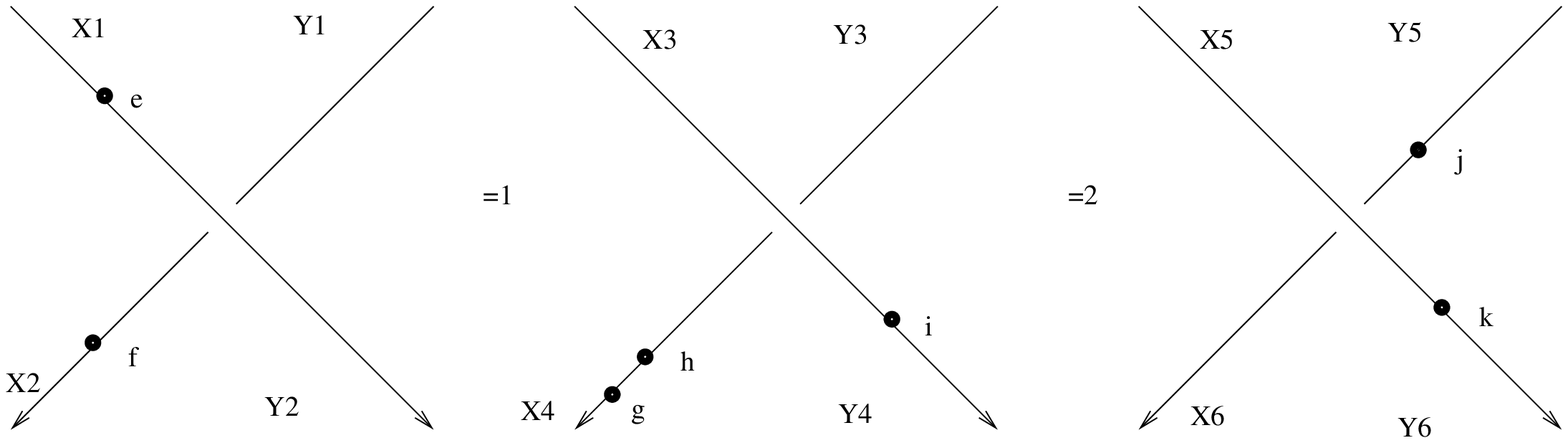}
\relabel{=1}{$=$}
\relabel{=2}{$=$}
\relabel{X1}{$\s{X}$}
\relabel{X2}{$\s{X}$}
\relabel{X3}{$\s{X}$}
\relabel{X4}{$\s{X}$}
\relabel{X5}{$\s{X}$}
\relabel{X6}{$\s{X}$}
\relabel{Y1}{$\s{\d(e)X}$}
\relabel{Y2}{$\s{\d(e)X}$}
\relabel{Y3}{$\s{\d(e)X}$}
\relabel{Y4}{$\s{\d(e)X}$}
\relabel{Y5}{$\s{\d(e)X}$}
\relabel{Y6}{$\s{\d(e)X}$}
\relabel{e}{$\s{e}$}
\relabel{f}{$\s{e^{-1}}$}
\relabel{g}{$\s{e^{-1}}$}
\relabel{h}{$\s{eX \t e^{-1}}$}
\relabel{i}{$\s{e}$}
\relabel{j}{$\s{e^{-1}}$}
\relabel{k}{$\s{e}$}
\relabel{=1}{$=$}
\relabel{=2}{$=$}
\endrelabelbox}
\caption{An identity.}
\label{id5}
\end{figure}

\begin{figure}
\centerline{\relabelbox 
\epsfysize 6cm
\epsfbox{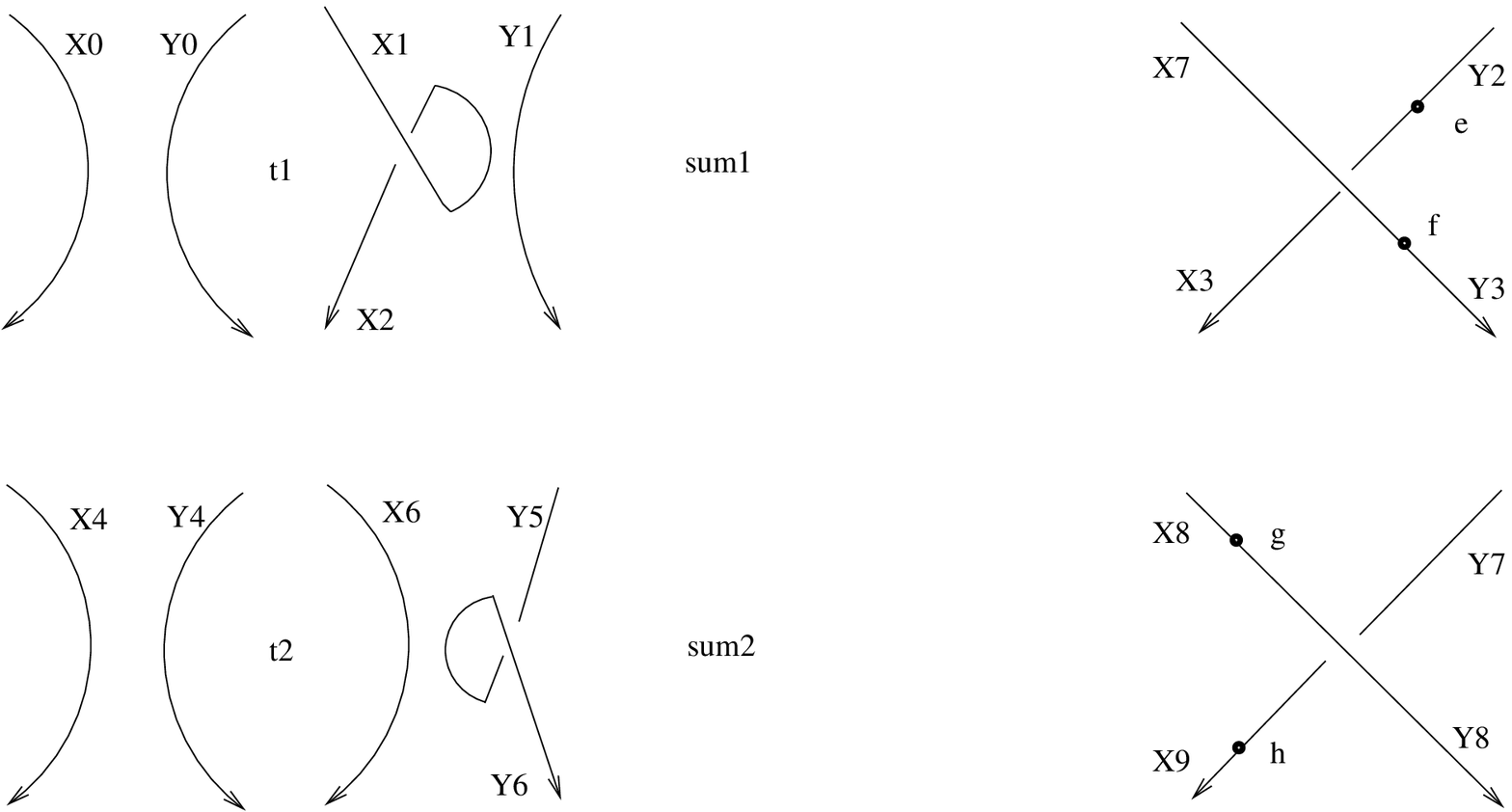}
\relabel{t1}{$\mapsto$}
\relabel{t2}{$\mapsto$}
\relabel{X0}{$\s{X}$}
\relabel{X1}{$\s{X}$}
\relabel{X2}{$\s{X}$}
\relabel{X3}{$\s{X}$}
\relabel{X4}{$\s{X}$}
\relabel{X7}{$\s{X}$}
\relabel{X6}{$\s{X}$}
\relabel{X8}{$\s{X}$}
\relabel{X9}{$\s{X}$}
\relabel{Y0}{$\s{Y}$}
\relabel{Y1}{$\s{Y}$}
\relabel{Y4}{$\s{Y}$}
\relabel{Y5}{$\s{Y}$}
\relabel{Y6}{$\s{Y}$}
\relabel{Y2}{$\s{\d(e)X}$}
\relabel{Y3}{$\s{\d(e)X}$}
\relabel{Y7}{$\s{\d(e)X}$}
\relabel{Y8}{$\s{\d(e)X}$}
\relabel{e}{$\s{e^{-1}}$}
\relabel{f}{$\s{e}$}
\relabel{g}{$\s{e}$}
\relabel{h}{$\s{ e^{-1}}$}
\relabel{sum1}{$\displaystyle{\mapsto \frac{1}{\# E}\sum_{e \in E} \delta(\d(e)X,Y)}$}
\relabel{sum2}{$\displaystyle{\mapsto \frac{1}{\# E}\sum_{e \in E} \delta(\d(e)X,Y)}$}
\endrelabelbox}
\caption{Invariance under Going Right Movie Move $13$.}
\label{MM13rightproof}
\end{figure}

The movie move $13$ is presented in figure \ref{MM13}. It should be considered
in both directions and considering also mirror images and opposite, though
compatible, orientations.
To begin with, we draw attention to the identity of figure \ref{id5}. The
proof of invariance under Going Right Movie Move $13$ is a corollary of this
identity, and appears in figure \ref{MM13rightproof}.
 
\begin{figure}
\centerline{\relabelbox 
\epsfysize 6cm
\epsfbox{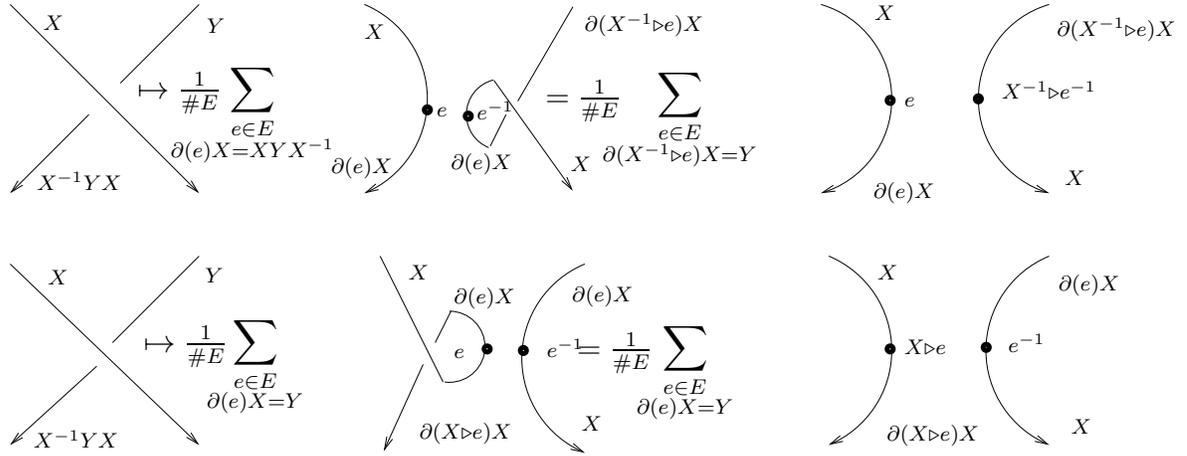}
\relabel{X0}{$\s{X}$}
\relabel{Y0}{$\s{Y}$}
\relabel{Z0}{$\s{X^{-1}YX}$}
\relabel{X1}{$\s{X}$}
\relabel{X4}{$\s{\d(e)X}$}
\relabel{Y1}{$\s{\d(X^{-1}\t e) X}$}
\relabel{X2}{$\s{X}$}
\relabel{Z}{$\s{\d(e)X}$}
\relabel{X3}{$\s{X}$}
\relabel{X5}{$\s{\d(e)X}$}
\relabel{X8}{$\s{\d(X^{-1}\t e) X}$}
\relabel{X6}{$\s{X}$}
\relabel{e}{$\s{e}$}
\relabel{f}{$\s{e^{-1}}$}
\relabel{h}{$\s{e}$}
\relabel{g}{$\s{X^{-1}\t e^{-1}}$}
\relabel{t1}{$\mapsto \frac{1}{\# E}$}
\relabel{sum1}{$\s{\displaystyle{\sum_{\substack{e \in E\\
          \d(e)X=XYX^{-1}}}}}$}
\relabel{t2}{$= \frac{1}{\# E}$}
\relabel{sum2}{$\s{\displaystyle{\sum_{\substack{e \in E\\
          \d(X^{-1}\t e)X=Y}}}}$}
\relabel{Y2}{$\s{X}$}
\relabel{Y3}{$\s{Y}$}
\relabel{Y4}{$\s{X^{-1}YX}$}
\relabel{Y5}{$\s{X}$}
\relabel{Z1}{$\s{\d(e)X}$}
\relabel{Y6}{$\s{\d(X\t e)X}$}
\relabel{Y7}{$\s{\d( e) X}$}
\relabel{Y8}{$\s{X}$}
\relabel{Y9}{$\s{X}$}
\relabel{Y10}{$\s{\d(e)X}$}
\relabel{Y11}{$\s{\d(X \t e) X}$}
\relabel{Y12}{$\s{X}$}
\relabel{i}{$\s{e}$}
\relabel{j}{$\s{e^{-1}}$}
\relabel{k}{$\s{X \t e}$}
\relabel{l}{$\s{e^{-1}}$}
\relabel{t3}{$\mapsto \frac{1}{\# E}$}
\relabel{sum3}{$\s{\displaystyle{\sum_{\substack{e \in E\\
          \d(e)X=Y}}}}$}
\relabel{t4}{$= \frac{1}{\# E}$}
\relabel{sum4}{$\s{\displaystyle{\sum_{\substack{e \in E\\
          \d( e)X=Y}}}}$}
\endrelabelbox}
\caption{Invariance under Going Left Movie Move $13$.}
\label{MM13left}
\end{figure}

The invariance under Going Left Movie Move $13$ appears in figure
\ref{MM13left}.

\subsubsection{Movie Move 14}
\begin{figure}
\centerline{\relabelbox
\epsfysize 3cm
\epsfbox{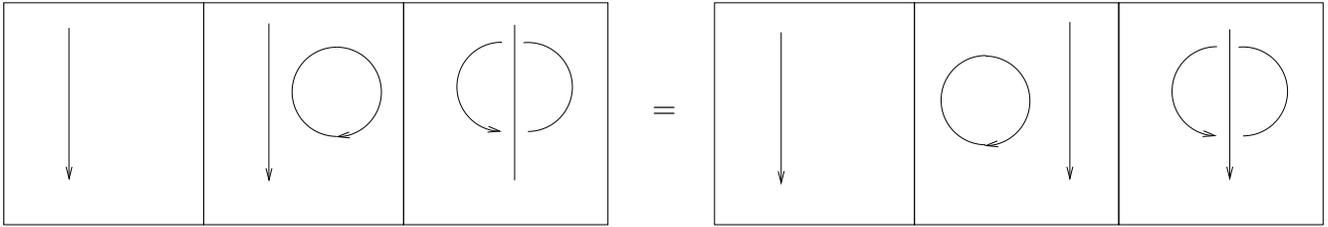}
\relabel{=}{$=$}
\endrelabelbox }
\caption{\label{MM14_good} A version of Movie Move $14$.}
\end{figure}

A version of the Movie Move 14 appears in figure \ref{MM14_good}.
This is a reversible movie, so we need to treat the Going Right and Going
Left cases separately. As usual, we need to consider all the possible changes in the
orientations  of the edges as well as on the  crossing information.

The going right case  is the simplest, since we can suppose that no vertices get in the way. We leave the proof to the reader (c.f. section \ref{trivial}).

\begin{figure}
\centerline{\relabelbox
\epsfysize 7cm
\epsfbox{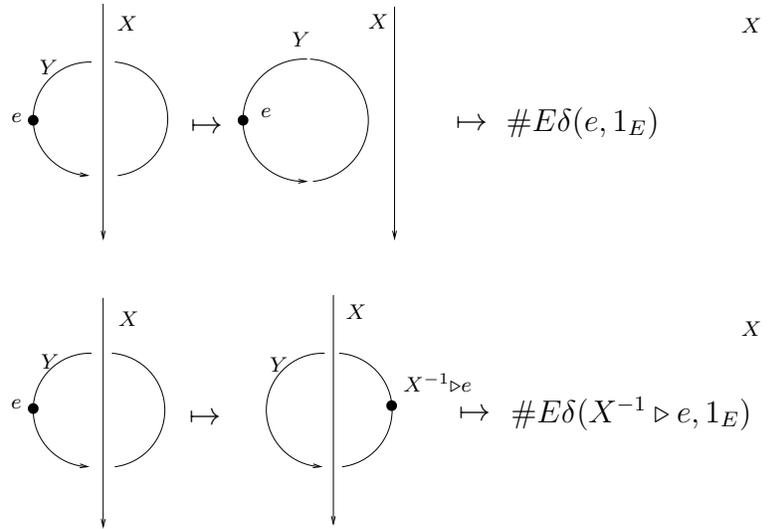}
\relabel{t1}{$\mapsto$}
\relabel{t2}{$\mapsto$}
\relabel{t3}{$\mapsto$}
\relabel{t4}{$\mapsto$}
\relabel{A}{$\s{X}$}
\relabel{C}{$\s{X}$}
\relabel{E}{$\s{X}$}
\relabel{F}{$\s{X}$}
\relabel{H}{$\s{X}$}
\relabel{J}{$\s{X}$}
\relabel{B}{$\s{Y}$}
\relabel{D}{$\s{Y}$}
\relabel{G}{$\s{Y}$}
\relabel{I}{$\s{Y}$}
\relabel{a}{$\s{e}$}
\relabel{b}{$\s{e}$}
\relabel{c}{$\s{e}$}
\relabel{d}{$\s{X^{-1}\t e}$}
\relabel{D1}{$\# E\delta(e,1_E)$}
\relabel{D2}{$\# E \delta(X^{-1}\t e,1_E)$}
\endrelabelbox }
\caption{\label{MM14_goodproof1} Invariance under Going Left Movie Move $14$, first case.}
\end{figure}

\begin{figure}
\centerline{\relabelbox
\epsfysize 7cm
\epsfbox{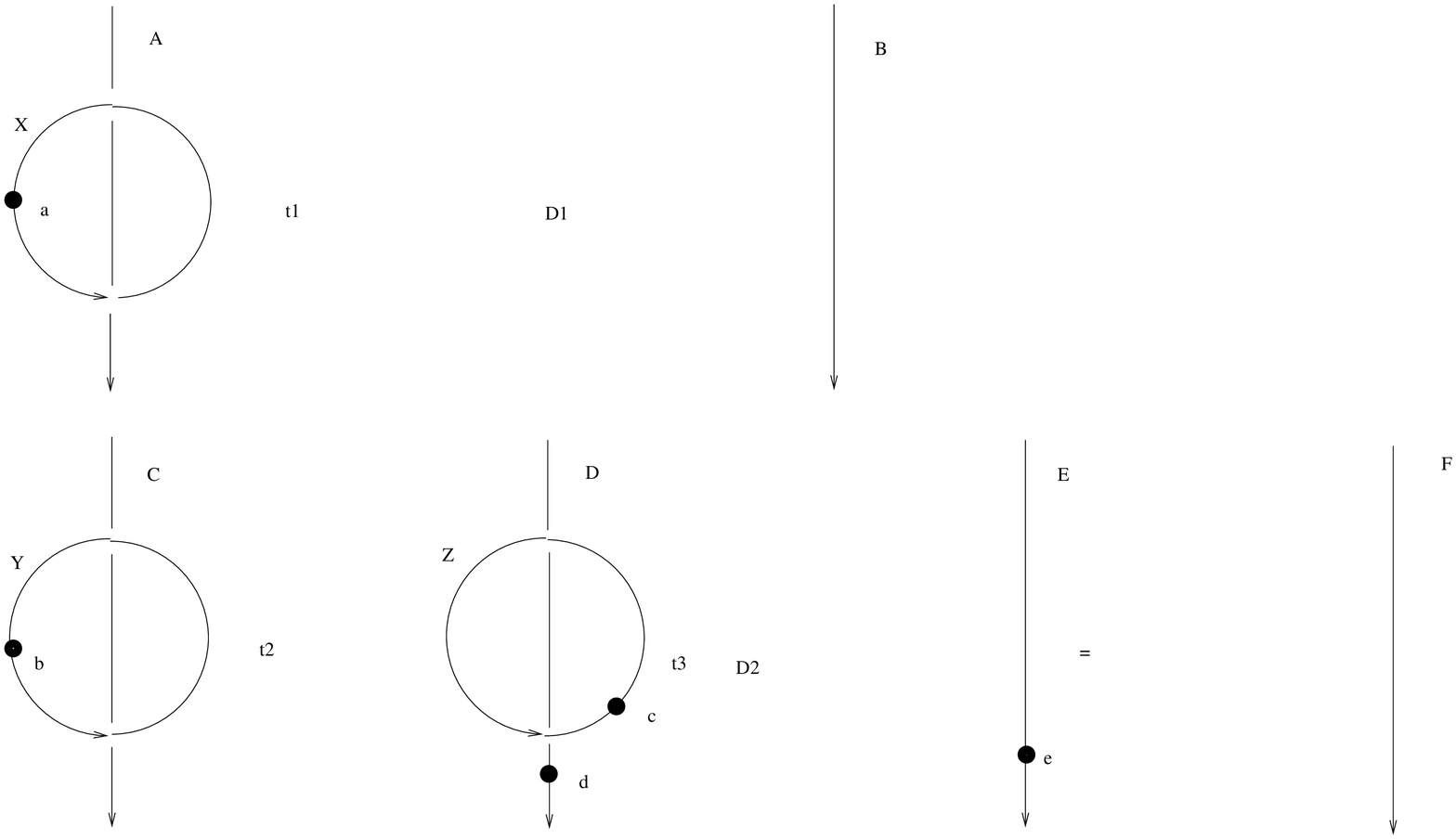}
\relabel{t1}{$\mapsto$}
\relabel{t2}{$\mapsto$}
\relabel{t3}{$\mapsto$}
\relabel{=}{$=\#E \delta(e,1_E)$}
\relabel{A}{$\s{X}$}
\relabel{B}{$\s{X}$}
\relabel{C}{$\s{X}$}
\relabel{D}{$\s{X}$}
\relabel{E}{$\s{X}$}
\relabel{F}{$\s{X}$}
\relabel{X}{$\s{Y}$}
\relabel{Y}{$\s{Y}$}
\relabel{Z}{$\s{Y}$}
\relabel{a}{$\s{e}$}
\relabel{b}{$\s{e}$}
\relabel{c}{$\s{e}$}
\relabel{d}{$\s{eXYX^{-1}\t e^{-1}}$}
\relabel{e}{$\s{eXYX^{-1}e^{-1}}$}
\relabel{D1}{$\#E \delta(e,1_E)$}
\relabel{D2}{$\#E \delta( e,1_E) $}
\endrelabelbox }
\caption{\label{MM14_goodproof2} Invariance under Going Left Movie Move $14$,
  second  case.}
\end{figure}

For the going left case, we can suppose that we initially have a unique vertex
in the circle. The invariance under the two different types of Going Left Movie
Move $14$ appears in figures \ref{MM14_goodproof1} and
\ref{MM14_goodproof2}. The cases with different orientations are analogous.

\subsubsection{Movie Move 15}
\begin{figure}
\centerline{\relabelbox 
\epsfysize 2cm
\epsfbox{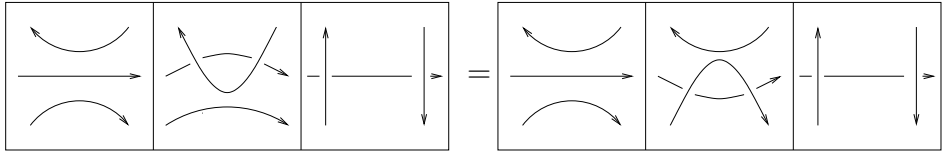}
\relabel {=}{$=$}
\endrelabelbox }
\caption{\label{MM14} A version of Movie Move $15$.}
\end{figure}

A version of the oriented Movie Move $15$ appears in figure \ref{MM14}. There
are several other versions obtained by reversing the orientations of the
strands and considering  mirror images. We say a movie move $15$ is of the
first kind if it is obtained from figure \ref{MM14} by (possibly) changing
the orientation of the strands, and it is of the second kind if, similarly,  it is obtained from the mirror image of figure \ref{MM14}. 

\begin{figure}
\centerline{\relabelbox
\epsfysize 3cm
\epsfbox{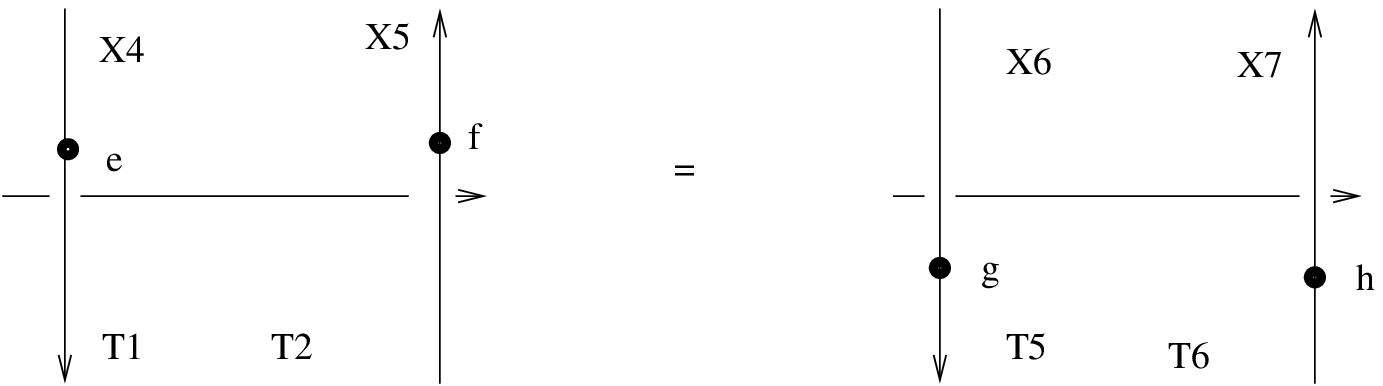}
\relabel{X4}{$\s{X}$}
\relabel{X5}{$\s{X}$}
\relabel{X6}{$\s{X}$}
\relabel{X7}{$\s{X}$}
\relabel{T1}{$\s{\d(e)X}$}
\relabel{T2}{$\s{\d(e)X}$}
\relabel{T5}{$\s{\d(e)X}$}
\relabel{T6}{$\s{\d(e)X}$}
\relabel{e}{$\s{e}$}
\relabel{f}{$\s{e^{-1}}$}
\relabel{g}{$\s{e}$}
\relabel{h}{$\s{e^{-1}}$}
\relabel{=}{$=$}
\endrelabelbox}
\caption{An identity.}
\label{id7}
\end{figure}

\begin{figure}
\centerline{\relabelbox 
\epsfysize 6cm
\epsfbox{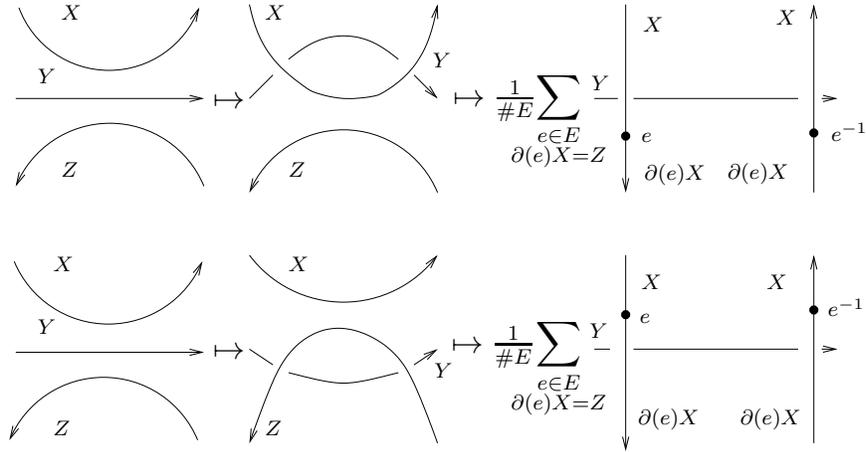}
\relabel{X4}{$\s{X}$}
\relabel{X5}{$\s{X}$}
\relabel{X6}{$\s{X}$}
\relabel{X7}{$\s{X}$}
\relabel{T1}{$\s{\d(e)X}$}
\relabel{T2}{$\s{\d(e)X}$}
\relabel{T3}{$\s{\d(e)X}$}
\relabel{T4}{$\s{\d(e)X}$}
\relabel{e}{$\s{e}$}
\relabel{f}{$\s{e^{-1}}$}
\relabel{g}{$\s{e}$}
\relabel{h}{$\s{e^{-1}}$}
\relabel{t1}{$\mapsto$}
\relabel{t2}{$\mapsto$}  
\relabel{t3}{$\mapsto\frac{1}{\# E}$}
\relabel{t4}{$\mapsto\frac{1}{\# E}$}
\relabel{X0}{$\s{X}$}
\relabel{X2}{$\s{X}$}
\relabel{X1}{$\s{X}$}
\relabel{X3}{$\s{X}$}
\relabel{Y0}{$\s{Y}$}
\relabel{Y1}{$\s{Y}$}
\relabel{Y2}{$\s{Y}$}
\relabel{Y3}{$\s{Y}$}
\relabel{Y4}{$\s{Y}$}
\relabel{Y5}{$\s{Y}$}
\relabel{Z0}{$\s{Z}$}
\relabel{Z1}{$\s{Z}$}
\relabel{Z2}{$\s{Z}$}
\relabel{Z3}{$\s{Z}$}
\relabel{sum1}{$\displaystyle{\sum_{\substack{ e\in E\\ \d(e)X=Z }}}$}
\relabel{sum2}{$\displaystyle{\sum_{\substack{ e\in E\\ \d(e)X=Z }}}$}
\endrelabelbox}
\caption{Invariance under Going Right Movie Move $15$, first case.}
\label{MM14right}
\end{figure}

\begin{figure}
\centerline{\relabelbox 
\epsfysize 8cm
\epsfbox{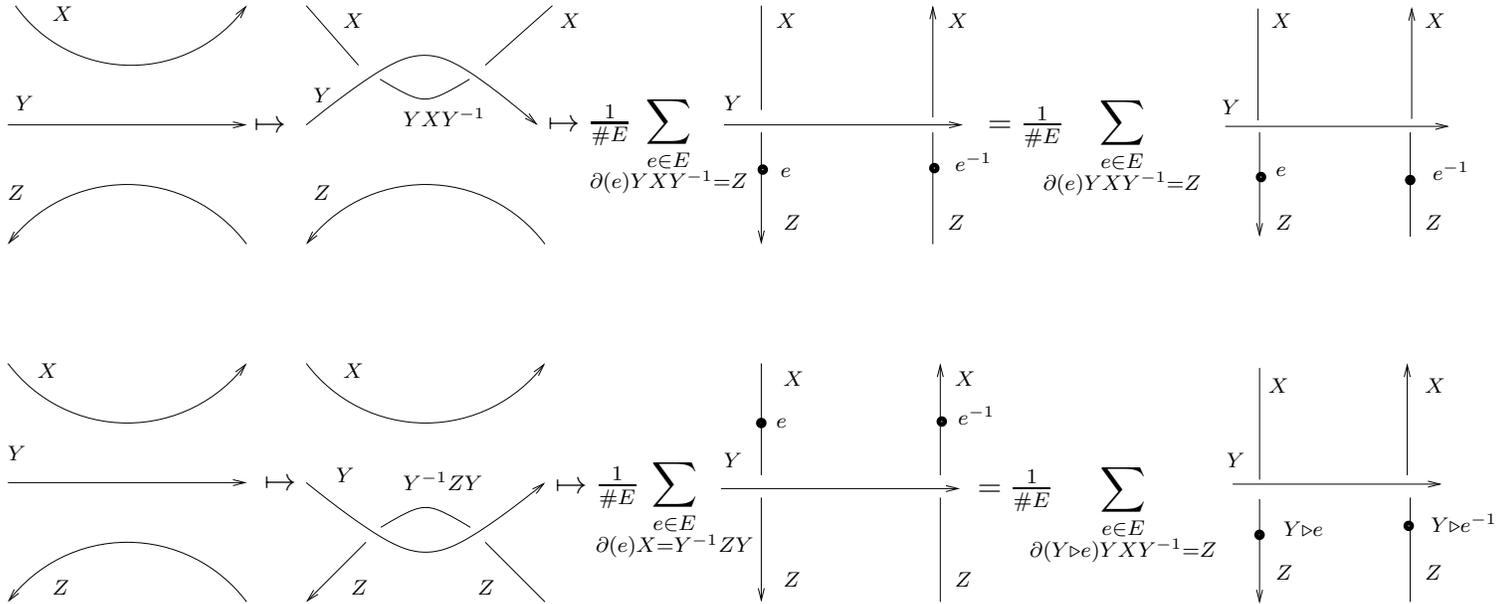}
\relabel{X0}{$\s{X}$}
\relabel{X1}{$\s{X}$}
\relabel{X2}{$\s{X}$}
\relabel{X3}{$\s{X}$}
\relabel{X4}{$\s{X}$}
\relabel{X5}{$\s{X}$}
\relabel{X6}{$\s{X}$}
\relabel{X7}{$\s{X}$}
\relabel{X8}{$\s{X}$}
\relabel{X9}{$\s{X}$}
\relabel{X10}{$\s{X}$}
\relabel{X11}{$\s{X}$}
\relabel{X12}{$\s{X}$}
\relabel{W1}{$\s{YXY^{-1}}$}
\relabel{Y0}{$\s{Y}$}
\relabel{Y1}{$\s{Y}$}
\relabel{Y2}{$\s{Y}$}
\relabel{Y3}{$\s{Y}$}
\relabel{Y4}{$\s{Y}$}
\relabel{Y5}{$\s{Y}$}
\relabel{Y6}{$\s{Y}$}
\relabel{Y7}{$\s{Y}$}
\relabel{W2}{$\s{Y^{-1}ZY}$}
\relabel{Z1}{$\s{Z}$}
\relabel{Z2}{$\s{Z}$}
\relabel{Z3}{$\s{Z}$}
\relabel{Z4}{$\s{Z}$}
\relabel{Z5}{$\s{Z}$}
\relabel{Z6}{$\s{Z}$}
\relabel{Z7}{$\s{Z}$}
\relabel{Z8}{$\s{Z}$}
\relabel{Z9}{$\s{Z}$}
\relabel{Z10}{$\s{Z}$}
\relabel{Z11}{$\s{Z}$}
\relabel{Z12}{$\s{Z}$}
\relabel{Z13}{$\s{Z}$}
\relabel{e}{$\s{e}$}
\relabel{f}{$\s{e^{-1}}$}
\relabel{g}{$\s{e}$}
\relabel{h}{$\s{e^{-1}}$}
\relabel{i}{$\s{e}$}
\relabel{j}{$\s{e^{-1}}$}
\relabel{k}{$\s{Y \t e}$}
\relabel{l}{$\s{Y \t e^{-1}}$}
\relabel{t1}{${\mapsto}$}
\relabel{t2}{${\mapsto}$}
\relabel{t3}{$  {\mapsto\frac{1}{\# E}}$}
\relabel{t4}{$  {\mapsto\frac{1}{\# E}}$}
\relabel{sum1}{$\displaystyle{\sum_{\substack{ e\in E\\ \d(e)YXY^{-1}=Z }}}$}
\relabel{sum2}{$\displaystyle{\sum_{\substack{ e\in E\\ \d(e)YXY^{-1}=Z }}}$}
\relabel{sum3}{$\displaystyle{\sum_{\substack{ e\in E\\ \d(e)X=Y^{-1}ZY }}}$}
\relabel{sum4}{$\displaystyle{\sum_{\substack{ e\in E\\ \d(Y \t e)YXY^{-1}=Z}}}$}
\relabel{=1}{$=\frac{1}{\# E}$}
\relabel{=2}{$=\frac{1}{\# E}$}
\endrelabelbox }
\caption{\label{MM14right2} Invariance under Going Right Movie Move $15$, second case.}
\end{figure}

Let us consider the going right case first.  If the movie move is of
the first kind, the invariance is a consequence of the identity in figure \ref{id7}, and the calculation appears in figure \ref{MM14right}.  The proof of invariance under the other kinds of going right Movie Move $15$ appears in figure \ref{MM14right2}.

\begin{figure}
\centerline{\relabelbox 
\epsfysize 4cm
\epsfbox{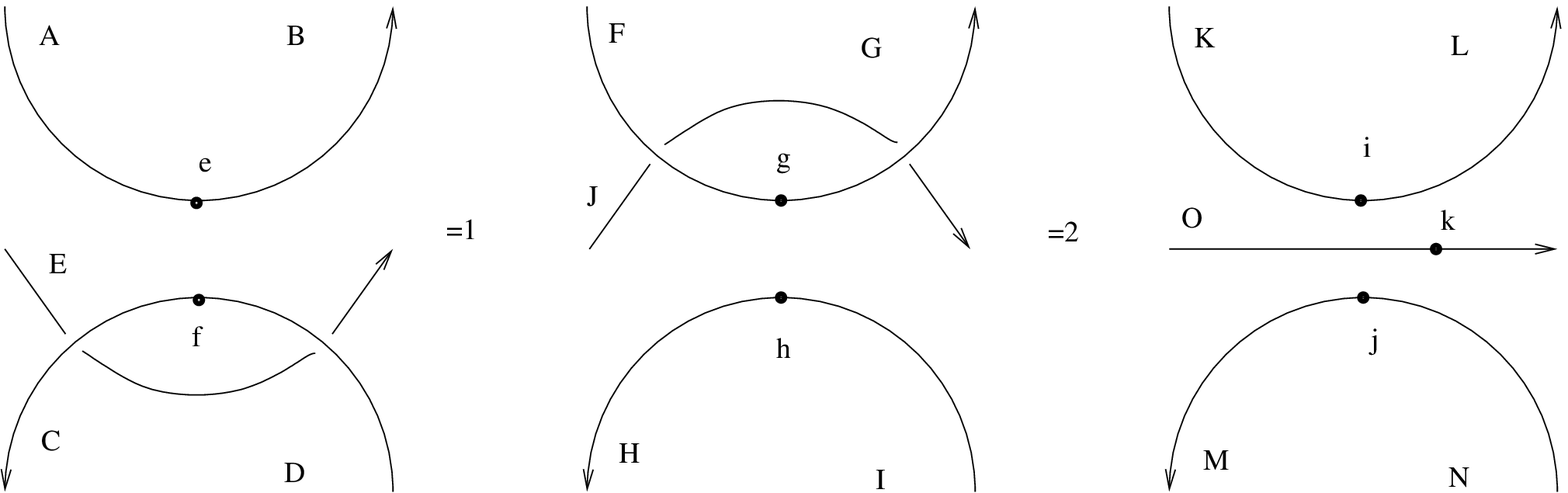}
\relabel{A}{$\s{X}$}
\relabel{C}{$\s{X}$}
\relabel{F}{$\s{X}$}
\relabel{H}{$\s{X}$}
\relabel{K}{$\s{X}$}
\relabel{M}{$\s{X}$}
\relabel{B}{$\s{\d(e)X}$}
\relabel{D}{$\s{\d(e)X}$}
\relabel{G}{$\s{\d(e)X}$}
\relabel{I}{$\s{\d(e)X}$}
\relabel{L}{$\s{\d(e)X}$}
\relabel{N}{$\s{\d(e)X}$}
\relabel{e}{$\s{e}$}
\relabel{e}{$\s{e}$}
\relabel{g}{$\s{e}$}
\relabel{i}{$\s{e}$}
\relabel{f}{$\s{e^{-1}}$}
\relabel{h}{$\s{e^{-1}}$}
\relabel{j}{$\s{e^{-1}}$}
\relabel{k}{$\s{e Z \t e^{-1}}$}
\relabel{E}{$\s{Z}$}
\relabel{J}{$\s{Z}$}
\relabel{O}{$\s{Z}$}
\relabel{=1}{$=$}
\relabel{=2}{$=$}
\endrelabelbox}
\caption{An identity.}
\label{id8}
\end{figure}

\begin{figure}
\centerline{\relabelbox 
\epsfysize 8cm
\epsfbox{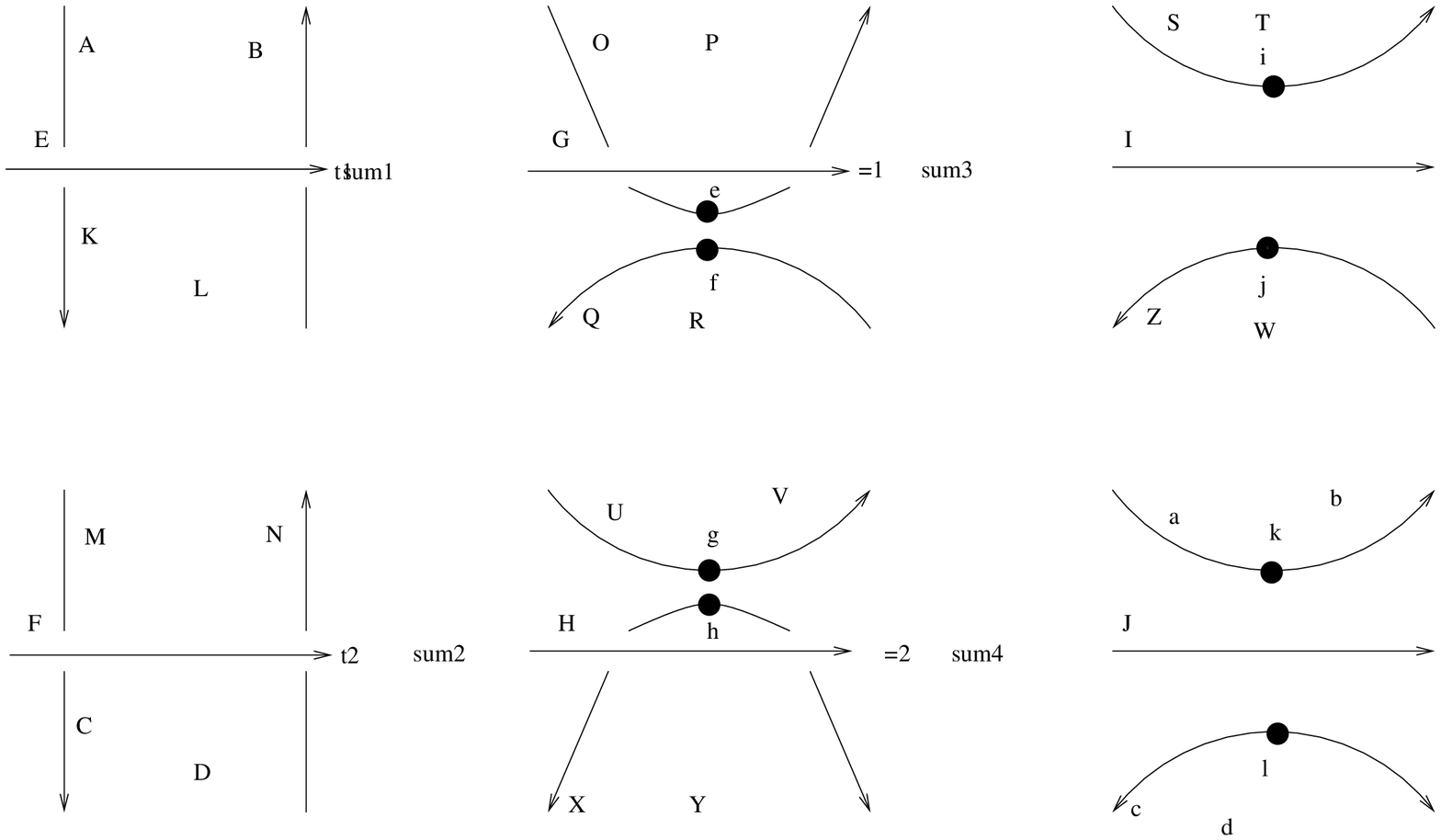}
\relabel{A}{$\s{X}$}
\relabel{O}{$\s{X}$}
\relabel{S}{$\s{X}$}
\relabel{M}{$\s{X}$}
\relabel{U}{$\s{X}$}
\relabel{a}{$\s{X}$}
\relabel{B}{$\s{Y}$}
\relabel{N}{$\s{Y}$}
\relabel{E}{$\s{Z}$}
\relabel{J}{$\s{Z}$}
\relabel{G}{$\s{Z}$}
\relabel{I}{$\s{Z}$}
\relabel{F}{$\s{Z}$}
\relabel{H}{$\s{Z}$}
\relabel{P}{$\s{\d(Z^{-1} \t e )X}$}
\relabel{T}{$\s{\d(Z^{-1} \t e )X}$}
\relabel{V}{$\s{\d(e)X}$}
\relabel{b}{$\s{\d(e)X}$}
\relabel{K}{$\s{ZXZ^{-1}}$}
\relabel{C}{$\s{ZXZ^{-1}}$}
\relabel{Q}{$\s{ZXZ^{-1}}$}
\relabel{X}{$\s{ZXZ^{-1}}$}
\relabel{Z}{$\s{ZXZ^{-1}}$}
\relabel{c}{$\s{ZXZ^{-1}}$}
\relabel{L}{$\s{ZYZ^{-1}}$}
\relabel{D}{$\s{ZYZ^{-1}}$}
\relabel{R}{$\s{\d(e)ZXZ^{-1}}$}
\relabel{W}{$\s{\d(e)ZXZ^{-1}}$}
\relabel{Y}{$\s{Z\d(e)XZ^{-1}}$}
\relabel{d}{$\s{\d(Z \t e)ZXZ^{-1}}$}
\relabel{=1}{$={\frac{1}{\# E}}$}
\relabel{=2}{$={\frac{1}{\# E}}$}
\relabel{e}{$\s{e}$}
\relabel{g}{$\s{e}$}               
\relabel{k}{$\s{e}$}
\relabel{f}{$\s{e^{-1}}$}
\relabel{h}{$\s{e^{-1}}$}
\relabel{j}{$\s{e^{-1}}$}
\relabel{i}{$\s{Z^{-1}\t e}$}
\relabel{l}{$\s{Z \t e^{-1}}$}
\relabel{t1}{$\mapsto{\frac{1}{\# E}}$}
\relabel{t2}{$\mapsto{\frac{1}{\# E}}$}
\relabel{sum1}{$\s{\displaystyle{\sum_{\substack{e\in E\\ \d(e)ZXZ^{-1}=ZYZ^{-1}   }} }}$}
\relabel{sum2}{$\s{\displaystyle{\sum_{\substack{e\in E\\ \d(e)X=Y   }} }}$}
\relabel{sum3}{$\s{\displaystyle{\sum_{\substack{e\in E\\ \d(Z^{-1} \t e)X=Y   }} }}$}
\relabel{sum4}{$\s{\displaystyle{\sum_{\substack{e\in E \\ \d(e)X=Y   }} }}$}
\endrelabelbox}
\caption{Invariance under the second kind of Going Left Movie Move 15.}
\label{MM14left2}
\end{figure}

Analogously, the invariance under the first kind of Going Left Movie Move $15$
is a consequence of the identity of figure \ref{id8}. The invariance under the
second kind of Going Left Movie Move $15$ appears in figure
\ref{MM14left2}. The cases with different orientation are dealt with similarly. This will finish the proof that $I(\S)$ is an isotopy invariant of knotted surfaces.

\begin{Theorem}
Let $\G=(G,E,\d,\t)$ be a finite crossed module. Let $\S$ be an oriented knotted surface. The evaluation $I^4_\G(\S)$ does not depend on the movie presentation of $\S$, and therefore $I^4_\G$ is an isotopy invariant of oriented knotted surfaces.
\end{Theorem}

\section{Discussion and Examples}\label{Discussion}

Let $\G=(G,E,\d,\t)$ be a finite categorical group. 
Suppose $E=\{1\}$ . For this particular case, in both $3$ and $4$ dimensions, the invariant $I=I_G$
is  simply the number of morphisms from  the
complement of the knot or the knotted surface into $G$, which tells us that
the invariant constructed in this article is non-trivial. The other extreme case with $G=\{1\}$ is less
interesting giving the trivial knot invariant in the $3$ dimensional case (consequence of the following proposition), and
only depending on the topology of the knotted surface in the $4$-dimensional
case. This is because, since $E=\ker \d$ is abelian, then the relations $R1$ to $R6$ imply that vertices can be moved freely, and considering that crossing information does not enter in the calculations, it is immediate that the invariant is defined only from  the handle decomposition of the knotted surface induced by the movie, thus not capturing  $4$-dimensional topology.

\begin{figure}
\centerline{\relabelbox 
\epsfysize 2cm
\epsfbox{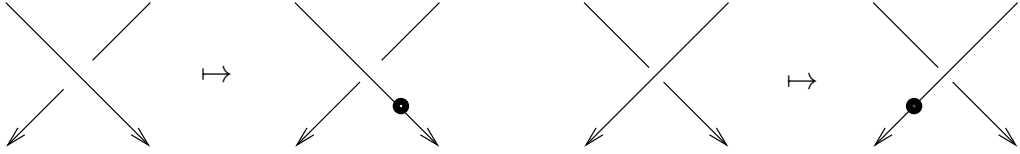}
\relabel{t1}{$\mapsto$}
\relabel{t2}{$\mapsto$}
\endrelabelbox}
\caption{Inserting new vertices in a knot.}
\label{K}
\end{figure}

\begin{figure}
\centerline{\relabelbox 
\epsfysize 4cm
\epsfbox{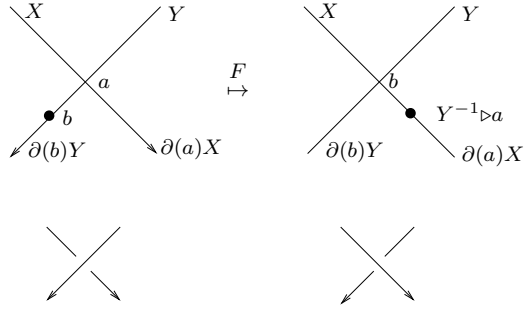}
\relabel{t}{$\substack{F\\ \mapsto}$}
\relabel{X}{$\s{X}$}
\relabel{A}{$\s{X}$}
\relabel{Y}{$\s{Y}$}
\relabel{B}{$\s{Y}$}
\relabel{Z}{$\s{\d(b)}Y$}
\relabel{W}{$\s{\d(a)X}$}
\relabel{C}{$\s{\d(b)Y}$}
\relabel{D}{$\s{\d(a)X}$}
\relabel{a}{$\s{a}$}
\relabel{b}{$\s{b}$}
\relabel{e}{$\s{b }$}
\relabel{g}{$\s{Y^{-1} \t a }$}
\endrelabelbox}
\caption{Map used to prove Proposition \ref{Refer5}}
\label{Kproof}
\end{figure}

A characteristic of the $3$-dimensional  case is:
\begin{Proposition}\label{Refer5}
Let $\G=(G,E,\d,\t)$ be a finite crossed module. In the $3$-dimensional case the invariant $I^3_\G$ is trivial whenever $G$ is abelian.\end{Proposition} 
\begin{Proof}
For a dotted link $K$, and one of its diagrams $D_K$,  let $K'$ be the dotted link obtained by inserting new vertices in $K$ as in figure \ref{K}. In particular
$I(K)=I(K')$ by \ref{vertex} (recall its results tell us that bivalent vertices can be absorbed by the rest of the diagram). Let $L$ be a link also constructed from $D_K$, but using different crossing information. We can also construct $L'$ by inserting new vertices in $L$, and we also have $I(L')=I(L)$. Let us prove that $I(L')=I(K')$, which implies that $I(K)=I(L)$. Let $C(D_{L'})$ and $C(D_{K'})$ be the sets of flat colourings of $D_{K'}$ and $D_{L'}$. Consider the map $F$ of figure \ref{Kproof}. This map is well defined since we have $\d(X\t e)=X \d(e) X^{-1}=\d(e)$.  It is immediate that  $F$ is a one-to-one map between the sets of flat colourings of $D_{L'}$ and  $D_{K'}$.
\end{Proof}

We prove in the rest of this section that, in the $4$ dimensional context, our invariant $I^4_\G$ is
stronger than the one $I_G$ made from a group $G$ alone, at least in a particular case. Let $\S\subset
S^4$ be a knotted sphere. It is well known, and easy to prove, that the first homology group of
the complement of $\S$ is always $\Z$. In particular if $G$ is an abelian
group then the  invariant $I_G$ is absolutely trivial for knotted spheres. We prove in the next subsection that there exists a crossed module $\G=(G,E,\d,\t)$ constructed using abelian groups such that the invariant $I^4_\G$  defined from  it 
detects the knottedness of the Spun Trefoil.

\subsection{Spun Trefoil}

\begin{figure}
\centerline{\relabelbox 
\epsfysize 5cm
\epsfbox{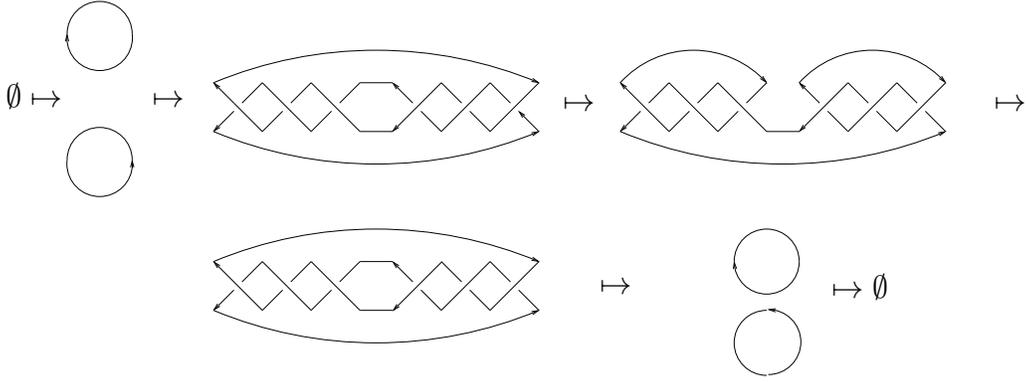}
\relabel{t1}{$\emptyset\mapsto$}
\relabel{t2}{$\mapsto$}
\relabel{t3}{$\mapsto$}
\relabel{t4}{$\mapsto$}
\relabel{t5}{$\mapsto$}
\relabel{t6}{$\mapsto\emptyset$}
\endrelabelbox}
\caption{A movie of the Spun Trefoil.}
\label{ST}
\end{figure}

\begin{figure}
\centerline{\relabelbox 
\epsfysize 6cm
\epsfbox{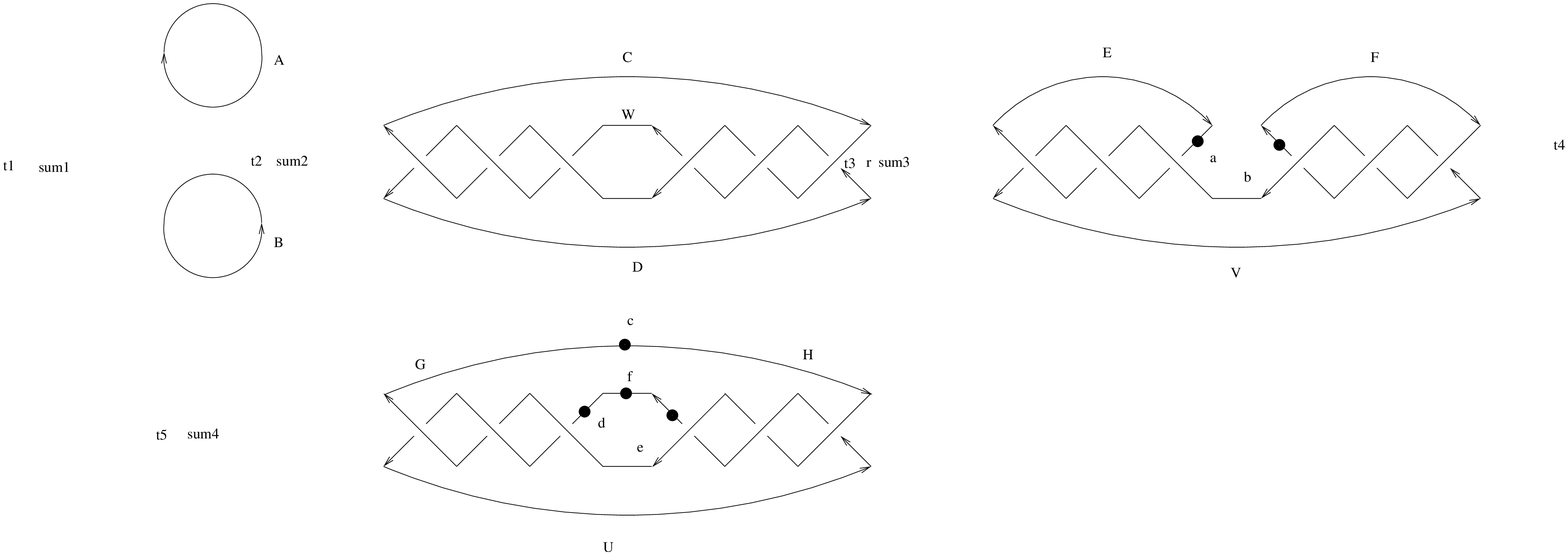}
\relabel{t1}{$1\mapsto$}
\relabel{t2}{$\mapsto$}
\relabel{t3}{$\mapsto$}
\relabel{t5}{$\mapsto$}
\relabel{sum1}{$\displaystyle{\sum_{X,Y \in G}}$}
\relabel{sum2}{$\displaystyle{\sum_{X,Y \in G}}$}
\relabel{r}{$\frac{1}{\#E}\ $}
\relabel{sum3}{$\displaystyle{\sum_{\substack{X,Y \in G\\ e \in
        E\\ \d(e)X=A}}}$}
\relabel{sum4}{$\frac{1}{(\#E)^2}\displaystyle{\sum_{\substack{X,Y \in G\\ e
        \in E\\ f \in \ker{\d} \\ d(e)X=A\ }}}$}
\relabel{A}{$\s{X}$}
\relabel{U}{$\s{Y}$}
\relabel{V}{$\s{Y}$}
\relabel{W}{$\s{A}$}
\relabel{B}{$\s{Y}$}
\relabel{C}{$\s{X}$}
\relabel{D}{$\s{Y}$}
\relabel{E}{$\s{X}$}
\relabel{F}{$\s{X}$}
\relabel{G}{$\s{X}$}
\relabel{H}{$\s{X}$}
\relabel{a}{$\s{e}$}
\relabel{b}{$\s{e^{-1}}$}
\relabel{c}{$\s{f}$}
\relabel{d}{$\s{e}$}
\relabel{f}{$\s{f^{-1}}$}
\relabel{e}{$\s{e^{-1}}$}
\endrelabelbox}
\caption{A calculation. Here we have $A=XYXYX^{-1}Y^{-1} X^{-1}$.}
\label{STinv}
\end{figure}

\begin{figure}
\centerline{\relabelbox 
\epsfysize 8cm
\epsfbox{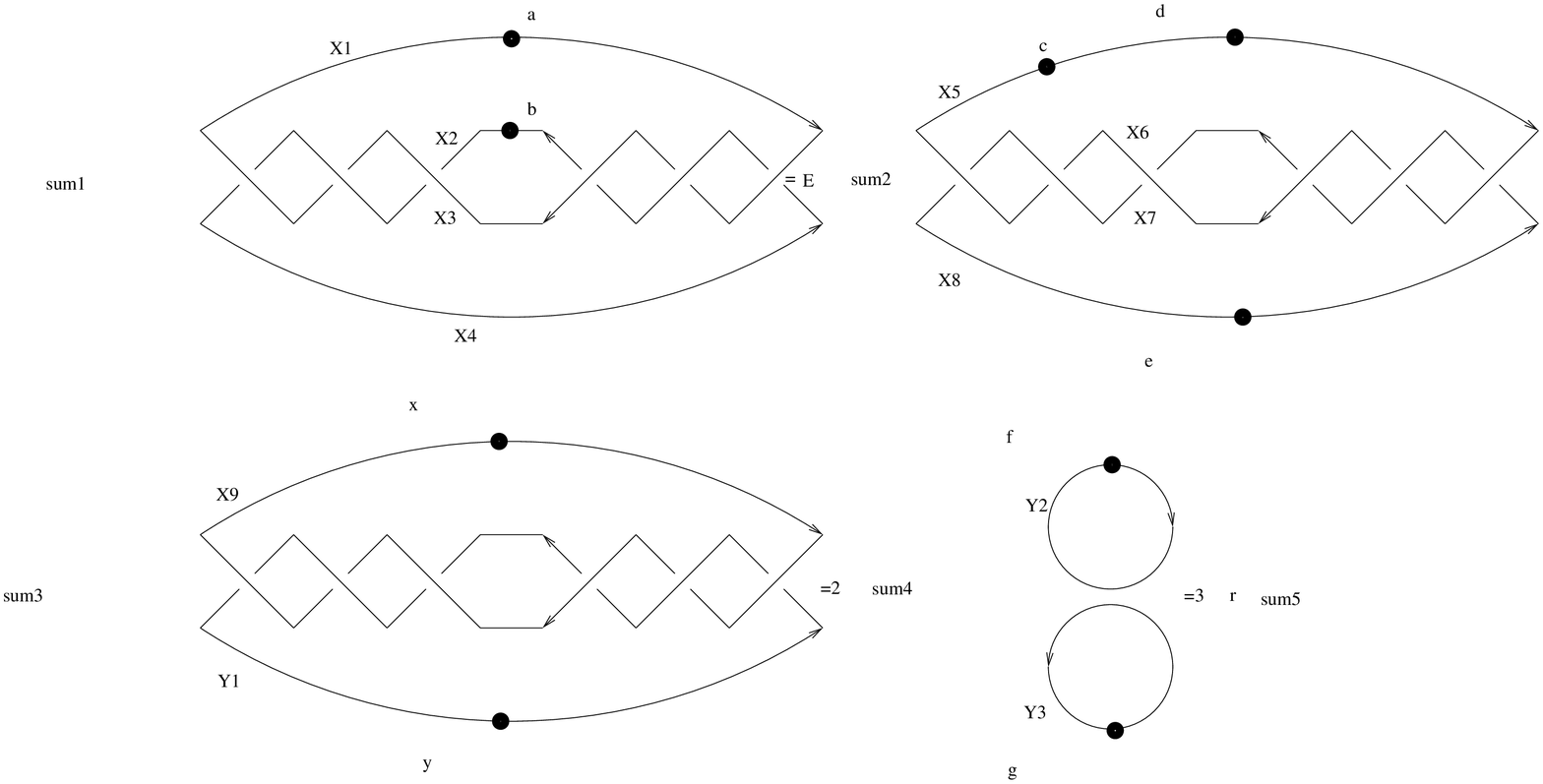}
\relabel{=}{$=$}
\relabel{=2}{$\mapsto$}
\relabel{=3}{$\mapsto$}
\relabel{E}{$\frac{1}{\#E}$}
\relabel{sum1}{$\frac{1}{\#E}\displaystyle{\sum_{\substack{X \in G\\ f \in E}}}$}
\relabel{sum2}{$\displaystyle{\sum_{\substack{X \in G\\ f \in E}}}$}
\relabel{sum3}{$=\frac{1}{\#E}\displaystyle{\sum_{\substack{X \in G\\ f \in E}}}$}
\relabel{sum4}{$\frac{1}{\#E}\displaystyle{\sum_{\substack{X \in G\\ f \in  E}}}$}
\relabel{X1}{$\s{X}$}
\relabel{X2}{$\s{X}$}
\relabel{X3}{$\s{X}$}
\relabel{X4}{$\s{X}$}
\relabel{X5}{$\s{X}$}
\relabel{X6}{$\s{X}$}
\relabel{X7}{$\s{X}$}
\relabel{X9}{$\s{X}$}
\relabel{Y1}{$\s{X}$}
\relabel{Y2}{$\s{X}$}
\relabel{Y3}{$\s{X}$}
\relabel{a}{$\s{f}$}
\relabel{b}{$\s{f^{-1}}$}
\relabel{c}{$\s{f}$}
\relabel{d}{$\s{X^{-1} \t f^{-1} X^{-2} \t f }$}
\relabel{e}{$\s{X^{-1} \t f^{-1} X^{-2} \t f X^{-3} \t f ^{-1}}$}
\relabel{x}{$\s{f X^{-1} \t f^{-1} X^{-2} \t f }$}
\relabel{y}{$\s{X^{-1}\left(f X^{-1} \t f^{-1} X^{-2} \t f  \right)^{-1} } $}
\relabel{f}{$\s{f X^{-1} \t f^{-1} X^{-2} \t f }$}
\relabel{g}{$\s{X^{-1}\left(f X^{-1} \t f^{-1} X^{-2} \t f  \right)^{-1} } $}
\relabel{r}{$(\#E)$}
\relabel{sum5}{$\displaystyle{\sum_{\substack{X \in G\\ f \in E\\f X^{-1} \t f^{-1} X^{-2} \t f=1_E}}} 1$}
\endrelabelbox}
\caption{Sequel of figure \ref{STinv} when both $G$ and $E$ are abelian, and $\d=1_E$.}
\label{STinvab}
\end{figure}

\begin{figure}
\centerline{\relabelbox 
\epsfysize 8cm
\epsfbox{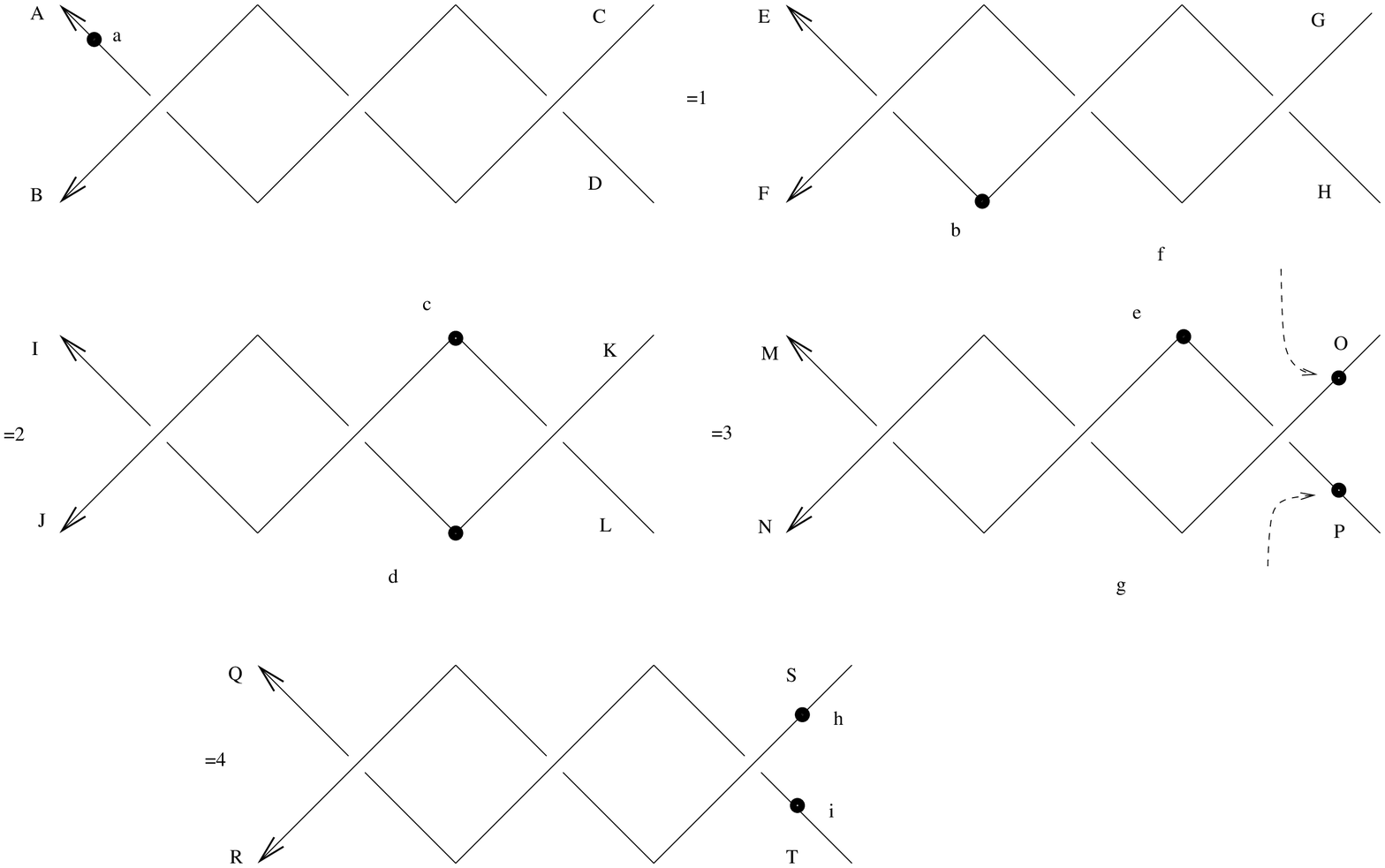}
\relabel{=1}{$=$}
\relabel{=2}{$=$}
\relabel{=3}{$=$}
\relabel{=4}{$=$}
\relabel{A}{$\s{X}$}
\relabel{B}{$\s{X}$}
\relabel{C}{$\s{X}$}
\relabel{D}{$\s{X}$}
\relabel{E}{$\s{X}$}
\relabel{F}{$\s{X}$}
\relabel{G}{$\s{X}$}
\relabel{H}{$\s{X}$}
\relabel{I}{$\s{X}$}
\relabel{J}{$\s{X}$}
\relabel{K}{$\s{X}$}
\relabel{L}{$\s{X}$}
\relabel{M}{$\s{X}$}
\relabel{N}{$\s{X}$}
\relabel{O}{$\s{X}$}
\relabel{P}{$\s{X}$}
\relabel{Q}{$\s{X}$}
\relabel{R}{$\s{X}$}
\relabel{S}{$\s{X}$}
\relabel{T}{$\s{X}$}
\relabel{a}{$\s{f^{-1}}$}
\relabel{b}{$\s{X^{-1}\t f^{-1}}$}
\relabel{c}{$\s{X^{-1}\t f^{-1}}$}
\relabel{d}{$\s{X^{-1}\t f^{-1} X^{-2}\t f }$}
\relabel{e}{$\s{X^{-1}\t f^{-1}}$}
\relabel{f}{$\s{X^{-1}\t f^{-1} X^{-2}\t f }$}
\relabel{g}{$\s{X^{-1} \t f^{-1} X^{-2} \t f X^{-3} \t f^{-1} X^{-2} \t f }$}
\relabel{h}{$\s{X^{-1} \t f^{-1} X^{-2} \t f }$}
\relabel{i}{$\s{X^{-1} \t f^{-1} X^{-3} \t f^{-1} X^{-2} \t f }$}
\endrelabelbox}
\caption{Proof of the second equality of the calculation in figure \ref{STinvab}. Note that all arcs are coloured with $X$.}
\label{calculation}
\end{figure}

A movie of the Spun Trefoil appears in figure \ref{ST}. We assign the obvious
movies to the second and the second last arrows. Observe that the choice of
these two movies is irrelevant for the final result, anyway. Let
$\G=(G,E,\d,\t)$ be a finite crossed module. Part of the  calculation of $I^4_\G(\S)$ for the Spun Trefoil appears in figure \ref{STinv}. 
Suppose  $G$ and $E$ are abelian and, moreover that $\d=1_G$ (notice that this implies that $E$ is abelian). We follow, for
this particular case, the calculation of figure \ref{STinv} in figure
\ref{STinvab}. The first equality in figure \ref{STinvab} is proved in figure \ref{calculation}. Therefore, if $G$ and $E$ are abelian, $\d=1_G$, and  $\S$ is the Spun Trefoil, we have:
$$I^4_\G(\S)=\#E\#\left \{(X,f): X \in G, f \in E: f X^{-1} \t f^{-1} X^{-2} \t f=1_E   \right \}.$$ 
On the other hand, it immediate that if $T$ is the trivial knotted surface we have (for any crossed module $\G$):
$$I^4_\G(T)=\#E \#G.$$
\begin{Remark}\label{Refer1}
Therefore, we will prove that the Spun Trefoil is knotted if we find a crossed module $\G=(G,E,\d,\t)$ made
from abelian groups, such that for some $X \in G$ there exists a solution of
$f X^{-1} \t f^{-1} X^{-2} \t f=1_E$ which is different from $f=1_E$.
\end{Remark}

\subsubsection{Explicit calculation for a particular crossed module}
Let $G$ and $E$ be groups with $E$ abelian. Suppose $G$ has a left action $\t$ on $E$ by automorphisms. If $\d=1_G$ then $\G=(G,E,\d,\t)$ is a crossed module.

Let $G$ be any group and $\k$ any field. Let $\k[G]$ be the free $\k$-vector
space on $G$. In particular $\k[G]$ is an abelian group and $G$ has an obvious
left action on $\k[G]$ by automorphisms. Therefore, we can construct a
crossed module if we are given a group $G$ and a field $\k$.

Consider $G=\Z_3=\{[0],[1],[2]\}$ and $\k=\Z_2$. We represent the elements of
$\Z_2[\Z_3]$ in the form $a.[0]+b.[1]+c.[2]$ where $a,b,c \in \Z_2$. Let
$f=1.[0]+1.[1]$ and $X=[1]$. We have $f- X^{-1} \t f +X^{-2} \t f=0$ (note that  we  have switched to additive notation, which is more appropriate for this example). For the
reasons pointed out in remark \ref{Refer1}, we thus conclude that, for this
case, the invariant $I^4_\G$ detects the knottedness of the Spun Trefoil.

The author is not aware  of  such a nice and simple example in the $3$-dimensional case. Note that Proposition \ref{Refer5} makes the task of finding simple examples  difficult.

\section*{Acknowledgements}
This work was done at the Mathematics Department of Instituto Superior T\'{e}cnico (Lisbon), with the financial support of FCT (Portugal), post-doc grant number SFRH/BPD/17552/2004. It is part of the research project POCTI/MAT/60352/2004 (''Quantum Topology''), also financed by FCT.

I would like to express my gratitude to Dr Roger Francis Picken and Dr Marco
Mackaay for many helpful discussions and their constant (and very encouraging) support from much before  this work had even started.  I would like
 to thank Dr Scott Carter for having been extremely helpful by answering numerous questions that appeared whilst I was preparing this article, and for
having commented on a preliminary version of it, encoranging me to calculate explicit examples.

This article shows the indirect, but strong,  influence of  Dr John W. Barrett, through the philosophy of working with handle decompositions instead of triangulations, a lesson, amongst many others, that I learnt from him whilst I was  his PhD student.

\end{document}